\documentclass[smallextended]{svjour3}       
\smartqed  
\usepackage{amsmath,amssymb,xcolor}
\usepackage{graphicx}
\usepackage{bbm}
\usepackage{float}
\usepackage{comment}
\usepackage{enumitem}
\usepackage{authblk}
\usepackage[font=footnotesize,labelfont=bf]{caption}
\usepackage{hyperref}
\usepackage{natbib}
\usepackage{multirow} 
\usepackage{footnote}
\usepackage{xspace}
\usepackage{makecell}

\newtheorem{condition}{Condition}

\newcommand\independent{\protect\mathpalette{\protect\independenT}{\perp}}
\def\independenT#1#2{\mathrel{\rlap{$#1#2$}\mkern2mu{#1#2}}}

\newcommand{\chr}{CHR\xspace}
\newcommand{\chrs}{CHRs\xspace}
\newcommand{\mchr}{MCHR\xspace}
\newcommand{\mchrs}{MCHRs\xspace}
\newcommand{\ohr}{OHR\xspace}
\newcommand{\ohrs}{OHRs\xspace}

\begin{document}

\title{The built-in selection bias of hazard ratios formalized}
\titlerunning{The built-in selection bias of hazard ratios formalized} 

\author{Richard A.J. Post         \and
 Edwin R. van den Heuvel \and
        Hein Putter 
       
}

\authorrunning{R.A.J. Post, E.R. van den Heuvel, and H. Putter} 

\institute{R.A.J. Post \at
              Department of Mathematics and Computer Science, Eindhoven University of Technology, The Netherlands \\
              \email{r.a.j.post@tue.nl}           
        \and 
        E.R. van den Heuvel \at
         Department of Mathematics and Computer Science, Eindhoven University of Technology, The Netherlands
        \and
           H. Putter \at
              Department of Biomedical Data Sciences, Leiden University Medical Center,  The Netherlands
}

\date{Received: date / Accepted: date}


	\maketitle
	
\begin{abstract}It is known that the hazard ratio lacks a useful causal interpretation. Even for data from a randomized controlled trial, the hazard ratio suffers from built-in selection bias as, over time, the individuals at risk in the exposed and unexposed are no longer exchangeable. In this work, we formalize how the observed hazard ratio evolves and deviates from the causal hazard ratio of interest in the presence of heterogeneity of the hazard of unexposed individuals (frailty) and heterogeneity in effect (individual modification). For the case of effect heterogeneity, we define the causal hazard ratio. We show that the observed hazard ratio equals the ratio of expectations of the latent variables (frailty and modifier) conditionally on survival in the world with and without exposure, respectively. Examples with gamma, inverse Gaussian and compound Poisson distributed frailty, and categorical (harming, beneficial or neutral) effect modifiers are presented for illustration. This set of examples shows that an observed hazard ratio with a particular value can arise for all values of the causal hazard ratio. Therefore, the hazard ratio can not be used as a measure of the causal effect without making untestable assumptions, stressing the importance of using more appropriate estimands such as contrasts of the survival probabilities.

\end{abstract}


\newpage
\section{Introduction}
When interested in time-to-event outcomes, ideally, one would like to know the hazard rates of an individual in the worlds with and without exposure. In practice, the focus is on the ratio of the expected hazard rates in both worlds. It is then standard practice to fit the observed hazard rates with a (time-invariant) Cox model \citep{Cox1972}. A decade ago, 
\citet{Hernan2010} raised awareness that hazard ratios estimated from a randomized controlled trial (RCT) are not suited for causal inference. Firstly, the average hazard ratio could be uninformative as there will typically be time-varying period-specific hazard ratios. More importantly, even if period-specific hazard ratios are estimated, these can vary solely due to the loss of randomization over time by conditioning on survivors. The exposure assignment and risk factors become dependent when conditioning on individuals that survived $t$, i.e.~survival time $T\geq t$, even if these risk factors are unrelated to the exposure \citep{Aalen2015}. As a result, effect measures based on hazards can suffer from non-collapsibility \citep{Daniel2020, Sjolander2016, Aalen2015, Martinussen2013}. 

In practice, the wrong estimand, a ratio of (partly) marginalized hazards, is defined, that by the non-collapsibility, deviates from the conditional (causal) hazard ratio. This contrast is referred to as the built-in selection bias of hazard ratios as the non-collapsibility results from conditioning on prior survival \citep{Hernan2010, Aalen2015, Sjolander2016, Stensrud2019, Young2020, Martinussen2020}. This bias should not be confused with confounding bias that is absent when using data from an RCT \citep{Didelez2021}. 
For exposure assignment $A$, and the potential survival time when the exposure is intervened on to $a$ denoted by $T^{a}$, the observed hazard ratio from an RCT satisfies
\begin{equation}\label{CH5hr:obs}
\frac{\lim_{h\rightarrow 0}h^{-1}\mathbb{P}\left(T \in [t,t+h) \mid T\geq t, A=a \right)}{\lim_{h\rightarrow 0}h^{-1}\mathbb{P}\left(T \in [t,t+h) \mid T\geq t, A=0 \right)}=\frac{\lim_{h\rightarrow 0}h^{-1}\mathbb{P}\left(T^{a} \in [t,t+h) \mid T^{a}\geq t \right)}{\lim_{h\rightarrow 0}h^{-1}\mathbb{P}\left(T^{0} \in [t,t+h) \mid T^{0}\geq t \right)}
\end{equation} \citep{DeNeve2020, Martinussen2020}. The observed hazard ratio thus equals the ratio of hazard rates at time $t$ for the potential outcomes of individuals from different populations; those for which $T^{a}\geq t$ and those for which $T^{0}\geq t$. As already indicated before, these populations will typically not be exchangeable in other risk factors, implying that an effect found cannot be (solely) assigned to the exposure. The effect does thus not reflect how the hazard rate of an individual is affected by exposure.  Only for cause-effect relations such that \eqref{CH5hr:obs} is time-invariant, the estimand can be interpreted as
\begin{equation*}
\frac{\log(\mathbb{P}(T^{a}\geq t))}{\log(\mathbb{P}(T^{0}\geq t))}.
\end{equation*}
It has therefore been recommended to use better interpretable estimands such as contrasts of quantiles, the restricted mean survival or survival probabilities of the potential outcomes respectively \citep{Hernan2010, Bartlett2020, Stensrud2019, Young2020}, or the probabilistic index derived from the latter \citep{DeNeve2020}. Alternatively, one can avoid interpretation issues by using accelerated failure time models \citep{Hernan2010, Hernan2005b} or additive hazard models \citep{Aalen2015, Martinussen2020}. 

Nevertheless, particularly in medical sciences, marginal hazard ratios are still commonly presented by practitioners. To boost the recommended paradigm shift, we would like to quantify the contrast between the causal hazard ratio and the observed hazard ratio, i.e.~the built-in selection bias. To do so, we start by presenting a general parameterization of cause-effect relations for time-to-event outcomes using a structural causal model. This parameterization clearly defines the causal effect of an exposure on an individual's hazard and, as such, allows us to define the causal hazard ratio. The quantitative examples in which the hazard under no exposure varies among individuals, i.e.~frailty, as presented in the literature \citep{Stensrud2017, Aalen2015, Balan2020} do fit in our framework, and we will formalize results for these examples. Additionally, we will extend these examples with causal effect heterogeneity, i.e.~the causal effect on the hazard rate might vary between individuals \citep{Stensrud2017}.  

\section{Notation}
In this work, probability distributions of factual and counterfactual outcomes are defined in terms of the potential outcome framework \citep{Rubin1974, Neyman1990}. Let $T_{i}$ and $A_{i}$ represent the (factual) stochastic outcome and exposure assignment level of individual $i$. Let $T_{i}^{a}$ equal the potential outcome of individual $i$ under an intervention of level $a$ (counterfactual when $A_{i} \neq a$). For those more familiar with the do-calculus,  $T^{a}$ is equivalent to $T \mid do(A=a)$ as e.g.~derived in \citet[Equation 40]{Pearl2009book} and \citet[Definition 8.6]{Bongers2021}. Throughout this work, we will assume causal consistency: if $A_{i}=a$, then $T_{i}^{a}=T_{i}^{A_{i}}=T_{i}$, implying that potential outcomes are independent of the exposure levels of other individuals (no interference). 

The hazard rate of $T^{a}$ can vary among the individuals in the population of interest. We will parameterize this heterogeneity for hazards of $T^{0}$ using a random variable $U_{0i}$ that represents the frailty of individual $i$ \citep{Balan2020}. The absolute effect of an exposure on the hazard, i.e.~the contrast between the hazard rates of $T_{i}^{0}$ and $T_{i}^{a}$ respectively, can depend on $U_{0i}$. However, there might also be (relative) effect heterogeneity that we parameterize using the random variable $U_{1i}$, giving rise to an individual-specific hazard ratio. The hazard of the potential outcome $T_{i}^{a}$ can be parameterized with a function that depends on $U_{0i}$, $U_{1i}$ and $a$. 
We focus on cause-effect relations that can be parameterized with a structural causal model (SCM)\footnote{A structural causal model as presented in this work is a union of the traditional SCM, for $a=A$, and the intervened SCMs for all possible $do(A=a)$ as e.g.~presented in \citet{Peters2017}.}, 
that consists of a joint probability distribution of $(N_{A}, U_{0}, U_{1}, N_{T})$ and a collection of structural assignments $(f_{A}, f_{\lambda})$ such that

\begin{center}
\fbox{%
	\parbox{0.9\linewidth}{%
		\begin{align}\label{CH5SCMsurv}
		A_{i}&:=f_{A}(N_{Ai})   \\\nonumber
		\lambda_{i}^{a}(t)&:= f_{\lambda}(t,U_{0i},U_{1i},a) \\ \nonumber
		T_{i}^{a}&:= \min \{t>0:e^{-\Lambda_{i}^{a}(t)}\leq N_{Ti}\},
		\end{align} 			where $\Lambda_{i}^{a}(t)=\int_{0}^{t}\lambda_{i}^{a}(s)ds$, $f_{\lambda}(t,U_{0i},U_{1i},0) \independent U_{1i}$, i.e.~when $a=0$ $f_{\lambda}$ does not contain $U_{1}$, and $N_{Ai},N_{Ti}\sim \text{Uni}[0,1]$, while $ U_{0i}, U_{1i} \independent N_{Ti}$ and $f_{A}$ is the inverse cumulative distribution function of $A$.}}\\ 
\end{center} 
\noindent Note that the data generating mechanism is described by this SCM as $T_{i}^{A_{i}}=T_{i}$. 	For the distribution of data observed from an RCT, we have that $N_{A}\independent U_{0}, U_{1}, N_{T}$, i.e.~no confounding, as a result of the randomization. Note that $\lambda_{i}^{a}(t)$ equals the hazard of the potential outcome of individual $i$ under exposure $a$, i.e.~$\lambda_{i}^{a}(t)=\lim_{h\rightarrow 0} h^{-1}\mathbb{P}\left(T_{i}^{a} \in [t,t+h) \mid T_{i}^{a} \geq t, U_{0i}, U_{1i} \right)$. In this parameterization, $U_{0}$ thus results in heterogeneity of the hazard under no exposure between individuals, and the presence of $U_{1}$ results in heterogeneity of the effect of the exposure on the hazard between individuals. The SCM could be re-parameterized by including more details, e.g.~measured risk factors, so that the remaining (conditional) unmeasured heterogeneity will be less variable. 

\begin{figure}[H]
	\centering
	\captionsetup{width=\textwidth}
	\includegraphics[scale=2.5]{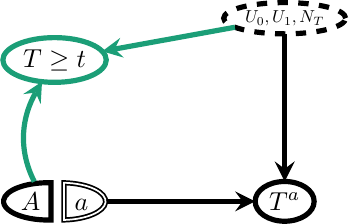}
	\caption{The SWIG (black) extended with the DAG (green) for SCM \eqref{CH5SCMsurv}, visualizing the causal relations between the exposure assignment $A$, the observed survival indicator $\{T\geq t\}$ and the potential survival time $T^{a}$. }\label{CH5CWIG}	
\end{figure}

Formulation of hazard rates of potential outcomes presented in the literature, e.g.~by \citet{Aalen2015} and by \citet{Stensrud2017}, naturally fit in this parameterization. However, the dependence of $T^{a}$ and $T^{0}$ beyond shared frailty is typically not specified. The causal relations in SCM \eqref{CH5SCMsurv} can be visualized in a graph. In Figure \ref{CH5CWIG}, one can observe this so-called Single-World Intervention Graph (SWIG) as proposed by \citet{Richardson2013} unified with the traditional causal directed acyclic graph (DAG) \citep{Pearl2009book}. 

It is insightful to realize that $T_{i}^{a}$ depends on $a$ and the random variables $U_{0i}$, $U_{1i}$, $N_{Ti}$ only, i.e.~ 
\begin{equation*}\label{CH5rand}
T_{i}^{a}:= \min \{t>0:e^{-\int_{0}^{t}f_{\lambda}(U_{0i},s,U_{1i},a)ds}\leq N_{Ti}\}=g(U_{0i},U_{1i},N_{Ti},a),
\end{equation*} for some function $g$, so by randomization
\begin{equation}\label{CH5inda}
T_{i}^{a} \independent A_{i}.
\end{equation}
However, it should be noted that
\begin{equation}\label{CH5indb}
T_{i}^{a} \not\independent A_{i} \mid T_{i} \geq  t, 
\end{equation}
as $T_{i}:=g(U_{0i}, U_{1i},N_{Ti}, A_{i})$, so that $A_{i} \mid T_{i}\geq t$ does inform on $(U_{0i},U_{1i}, N_{Ti})$ and thus on $T_{i}^{a}$. In the literature, \eqref{CH5indb} is often implicitly derived by recognizing that $\{T\geq t \}$ is a collider in Figure \ref{CH5CWIG}, and thus opens a back-door path between $A$ and $T^{a}$ \citep{Aalen2015, Sjolander2016}. The bias that results from conditioning on this collider is often referred to as the built-in selection bias of the hazard ratio. More general, in causal inference, selection bias is caused by conditioning on a variable $G$ that induces dependence between $Y^{a}$ and $A$, i.e.~$Y^{a} \independent A$ but $Y^{a} \not\independent A \mid G$, see e.g.~\citet[Chapter 8]{Hernan2019}. By causal consistency $Y^{a} \mid \{G, A=a\} \overset{d}{=} Y \mid G, A=a$ but $Y^{a} \mid G \overset{d}{\neq} Y^{a} \mid \{G, A=a\}$, so $Y^{a} \mid G \overset{d}{\neq} Y \mid G, A=a$ (as is the case in Equation \eqref{CH5indb} for $G=\{T\geq t\}$). However, the built-in selection bias is subtle as the dependence is a consequence of conditioning a potential outcome on information on its factual outcome, while it should be noted that (for an RCT),
\begin{equation}\label{CH5ind1}
T_{i}^{a} \independent A_{i} \mid T_{i}^{a}\geq  t, 
\end{equation} again as $T_{i}^{a}:=g(U_{0i}, U_{1i},N_{Ti}, a)$, $A_{i} \mid T^{a}_{i}\geq t$ does not inform on $(U_{0i},U_{1i}, N_{Ti})$. A more standard example of selection bias for survival analysis from an RCT is not accounting for informative censoring $(C=1)$ \citep{Hernan2004a, Howe2016}, in which case
\begin{equation*}
T^{a} \mid \{C=0 \}\overset{d}{\neq} T^{a} \mid \{C=0, A=a\}. 
\end{equation*} When interested in the causal effect, we should try to express the distribution of potential outcomes in terms of the observed distribution as we will do in the next section.

\section{The causal hazard ratio}\label{CH5sec:chr}
Parameterization of the cause-effect relations with SCM \eqref{CH5SCMsurv} allows us to define the causal hazard ratio.  If
\begin{equation*}
f_{\lambda}(t,U_{0i},U_{1i},a) = f_{0}(t,U_{0i})f_{1}(t,U_{1i},a) \text{ and } f_{1}(t,U_{1i},0)=1,
\end{equation*}
then the individual causal effect is captured by $f_{1}(t,U_{1i},a)$. The latter equals the ratio of the hazard of an individual's potential outcome when exposed and when non-exposed at time $t$, i.e.
$ \frac{\lambda_{i}^{a}(t)}{\lambda_{i}^{0}(t)}$. In the case of homogeneous effects, $f_{1}(t,U_{1i},a)=f_{1}(t,a)$ is equal for all individuals and would thus be the estimand of interest. In the case of heterogeneity of effects, $f_{1}(t,U_{1i},a)$ is the individual multiplicative causal effect. From a public health perspective, the ratio of the expected hazard in the world where everyone is exposed to $a$ and the expected hazard in the world where all individuals are unexposed is of interest. This causal hazard ratio (\chr) of interest can be obtained as the ratio of the marginalized (over $U_{0}$ and $U_{1}$) conditional hazard rates in both worlds as presented in Definition \ref{CH5def:CHR}.  \begin{definition}{\textbf{Causal hazard ratio}}\label{CH5def:CHR}
	The causal hazard ratio (\chr) for cause-effect relations that can be parameterized with SCM \ref{CH5SCMsurv} equals
	\begin{align}
 \frac{\mathbb{E}\left[\lambda_{i}^{a}(t)\right]}{\mathbb{E}\left[\lambda_{i}^{0}(t)\right]}&=\frac{\mathbb{E}\left[\lim_{h\rightarrow 0}h^{-1}\mathbb{P}\left(T^{a} \in [t,t+h) \mid T^{a}\geq t, U_{0}, U_{1}  \right)\right]}{\mathbb{E}\left[\lim_{h\rightarrow 0}h^{-1}\mathbb{P}\left(T^{0} \in [t,t+h) \mid T^{0}\geq t, U_{0} \right)\right]} \\ \nonumber 
 &= 
		\frac{\int\lim_{h\rightarrow 0}h^{-1}\mathbb{P}\left(T^{a} \in [t,t+h) \mid T^{a}\geq t, U_{0}, U_{1} \right)dF_{U_{0}, U_{1}}}{\int\lim_{h\rightarrow 0}h^{-1}\mathbb{P}\left(T^{0} \in [t,t+h) \mid  T^{0} \geq t, U_{0}  \right)dF_{U_{0}}},
	\end{align} where we abbreviate the Lebesque-Stieltjes integral of a function $g$ with respect to probability law $F_{X}$, i.e.~$\int g(x) dF_{X}(x)$, as $\int g(X) dF_{X}$. 
\end{definition} \noindent Note that, when $U_{0} \not \independent U_{1}$, the \chr can differ from the expectation of an individual's multiplicative effect ($\mathbb{E}\left[\frac{\lambda_{i}^{a}(t)}{\lambda_{i}^{0}(t)}\right]$). For cause-effect relations of a binary exposure that can be parameterized with SCM \eqref{CH5SCMsurv} such that $
f_{\lambda}(t,U_{0i},U_{1i},a) = U_{0i}\lambda_{0}(t)(U_{1i})^{a}$
 and $\mathbb{E}[U_{0}]=1$, $\frac{\mathbb{E}\left[\lambda_{i}^{1}(t)\right]}{\mathbb{E}\left[\lambda_{i}^{0}(t)\right]}=\mathbb{E}[U_{0}U_{1}]=\mathbb{E}[U_{1}]+\text{cov}(U_{0},U_{1})$ while $\mathbb{E}\left[\frac{\lambda_{i}^{1}(t)}{\lambda_{i}^{0}(t)}\right]=\mathbb{E}[U_{1}]$.\footnote{Consider for example, 
\begin{equation*}
(U_{0},U_{1}) = \begin{cases}
(0.5, 1.5) & p=0.5\\
(1.5, 0.5) & p=0.5,
\end{cases}
\end{equation*} so that apart from  $\mathbb{E}[U_{0}]=1$ also $\mathbb{E}[U_{1}]=1$ and $U_{0}U_{1}=\tfrac{3}{4}$ for everyone. As a result of the exposure, the hazard of $50\%$ of the individuals would increase by $50\%$ ($0.5\lambda_{0}(t) \rightarrow 0.75\lambda_{0}(t)$) and the hazard of the other $50\%$ would decrease by $50\%$ ($1.5\lambda_{0}(t) \rightarrow 0.75\lambda_{0}(t)$). The expected hazard rate in the world without exposure equals $\lambda_{0}(t)$ while the (expected) hazard in the world with exposure equals $0.75\lambda_{0}(t)$. On average, living in a world with exposure is thus more beneficial.} 
 
When the parameterization of the cause-effect relations as SCM \eqref{CH5SCMsurv} would be known, the \chr can be expressed in terms of the distribution of the data generating mechanism as presented in Theorem \ref{CH5th2}. \begin{theorem}\label{CH5th2}
	If the cause-effect relations of interest can be parameterized with SCM \eqref{CH5SCMsurv}, and $N_{A} \independent U_{0}, U_{1}, N_{T}$ (no confounding) then 
	\begin{equation*}
	\frac{\mathbb{E}\left[\lambda_{i}^{a}(t)\right]}{\mathbb{E}\left[\lambda_{i}^{0}(t)\right]}=\frac{\int\lim_{h\rightarrow 0}h^{-1}\mathbb{P}\left(T \in [t,t+h) \mid T\geq t, U_{0}, U_{1},A=a \right)dF_{U_{0}, U_{1}}}{\int\lim_{h\rightarrow 0}h^{-1}\mathbb{P}\left(T \in [t,t+h) \mid  T \geq t, U_{0}, A=0  \right)dF_{U_{0}}}.
	\end{equation*} 
\end{theorem} \noindent Consider, for example, the commonly used frailty model where effect heterogeneity is absent, i.e.
\begin{equation*}
\lambda_{i}^{a}(t) = U_{0i}\lambda_{0}(t)f_{1}(t,a).
\end{equation*} The \chr equals the multiplicative effect that does not differ among individuals and equals $f_{1}(t,a)$. By applying Theorem \ref{CH5th2}, this \chr is indeed derived to equal
\begin{equation}\label{CH5HRHom}
\frac{\int\lim_{h\rightarrow 0} h^{-1}\mathbb{P}\left(T \in [t,t+h) \mid T\geq t, U_{0}, A=a  \right)dF_{U_{0}}}{\int\lim_{h\rightarrow 0} h^{-1}\mathbb{P}\left(T \in [t,t+h) \mid T\geq t, U_{0}, A=0 \right)dF_{U_{0}}}=\frac{f_{1}(t,a)\lambda_{0}(t)\mathbb{E}[U_{0}]}{\lambda_{0}(t)\mathbb{E}[U_{0}]}=f_{1}(t,a).
\end{equation} It is important to note that \eqref{CH5HRHom} deviates from the observed hazard ratio equal to
\begin{equation*}
\frac{\lim_{h\rightarrow 0} h^{-1}\mathbb{P}\left(T \in [t,t+h) \mid T\geq t, A=a  \right)}{\lim_{h\rightarrow 0} h^{-1}\mathbb{P}\left(T \in [t,t+h) \mid T\geq t, A=0 \right)}=\frac{f_{1}(t,a)\lambda_{0}(t)\mathbb{E}[U_{0}\mid T\geq t, A=a]}{\lambda_{0}(t)\mathbb{E}[U_{0}\mid T\geq t, A=0]},
\end{equation*} as we will elaborate on in Section \ref{CH5OHR}. 

In summary, it becomes clear that to derive the \chr from data, inference on the distribution of the latent frailty $U_{0}$ and effect modifier $U_{1}$ must be made. When their distributions are known, inference can be drawn from observed data, even when the parameters of the distributions are unknown. Software available to estimate frailty parameters are described by \citet{Balan2020}, and such methods could also be adapted to estimate the latent modifier distribution. However, in practice, the distributions of these latent variables are unknown. Even in the case without causal effect heterogeneity, it is impossible to distinguish the presence of frailty from a time-dependent causal effect \citep[Section 2.5]{Balan2020}. More precisely, different combinations of (varying) effect sizes and frailty distributions give rise to the same marginal distribution. The same holds for combinations that also involve effect modifiers. 
In the case of clustered survival data, at least theoretically, the shared frailty could be distinguished from violation of proportional hazards \citep{Balan2020}. Reasoning along the same lines, individual frailty and marginal time-varying effects could only be derived from effect heterogeneity in the case of recurrent events with stationary distributions.

\section{Marginal causal hazard ratio}\label{CH5OHR}

In Theorem \ref{CH5th2}, the actual \chr (see Definition \ref{CH5def:CHR}) has been expressed in terms of the distribution of the observed data, and we concluded that these are typically not identifiable. Practitioners often compute the hazard ratio from data, which expectation is referred to as the observed hazard ratio (\ohr) and, in the absence of informative censoring, equals 
\begin{equation}\label{CH5eq:OH}
\text{OHR}(t)= \frac{\lim_{h\rightarrow 0}h^{-1}\mathbb{P}\left(T \in [t,t+h) \mid T\geq t, A=a \right)}{\lim_{h\rightarrow 0}h^{-1}\mathbb{P}\left(T \in [t,t+h) \mid T\geq t, A=0 \right)}.
\end{equation} To compare the \ohr to the \chr that quantifies the causal effect of interest, the \ohr should be expressed in terms of potential outcomes. For data from an RCT, by independence \eqref{CH5ind1} and causal consistency, the \ohr equals the marginal causal hazard ratio (\mchr), i.e.~
\begin{equation}\label{CH5nhr}
\text{MCHR}(t)=\frac{\lim_{h\rightarrow 0}h^{-1}\mathbb{P}\left(T^{a} \in [t,t+h) \mid T^{a}\geq t \right)}{\lim_{h\rightarrow 0}h^{-1}\mathbb{P}\left(T^{0} \in [t,t+h) \mid T^{0}\geq t \right)}. \end{equation} This \mchr should not be confused with the `marginal causal hazard ratio'  defined by \citet{Martinussen2020} that equals 
\begin{equation}
\frac{\lim_{h\rightarrow 0}h^{-1}\mathbb{P}\left(T^{a} \in [t,t+h) \mid T^{a}\geq t, T^{0}\geq t \right)}{\lim_{h\rightarrow 0}h^{-1}\mathbb{P}\left(T^{0} \in [t,t+h) \mid T^{0}\geq t, T^{a}\geq t \right)},
\end{equation} and is not considered in this work. 

We will study how the \mchr (and thus the \ohr from an RCT) differs from the \chr over time. By the law of total probability, \eqref{CH5nhr} equals
\begin{equation}
\frac{\int\lim_{h\rightarrow 0}h^{-1}\mathbb{P}\left(T^{a} \in [t,t+h) \mid T^{a}\geq t, U_{0}, U_{1} \right)dF_{U_{0}, U_{1}\mid T^{a}\geq t}}{\int\lim_{h\rightarrow 0}h^{-1}\mathbb{P}\left(T^{0} \in [t,t+h) \mid T^{0}\geq t, U_{0}\right)dF_{U_{0}\mid T^{0}\geq t}}. 
\end{equation} As the integration in the result of Theorem \ref{CH5th2} is with respect to the distribution of $U_{0}$ and $U_{1}$, in the entire population, instead of those individuals for which $T^{a}\geq t$ or $T^{0}\geq t$, the \mchr deviates from the \chr, i.e.~the built-in selection bias of the hazard \citep{Hernan2010, Aalen2015, Stensrud2019}. 

The problem induced for estimation of the \chr thus results from inference on the wrong estimand; the combined effect of the exposure of interest and the difference in latent frailty (and effect modification) distribution, is computed. This section will formalize how \eqref{CH5nhr} deviates from the \chr. Let us again focus on cause-effect relations that can be parameterized with SCM \eqref{CH5SCMsurv} in which
\begin{equation*}
\lambda_{i}^{a}(t) = f_{0}(t,U_{0i})f_{1}(t,U_{1i},a), 
\end{equation*} where the \chr equals $\mathbb{E}[f_{1}(t,U_{1},a)]$. These would be the cause-effect relations for which the \chr is an appropriate causal effect measure. Moreover, we focus on hazard functions that satisfy Condition \ref{CH5regcond} and do thus not have an infinite discontinuity.

\begin{condition}{\textbf{Hazard without infinite discontinuity}}\label{CH5regcond}
	$$\exists \tilde{{}h}>0 \text{ so that }\forall h^{*} \in (0,\tilde{{}h}):  \mathbb{E}\left[f_{0}(t+h^{*},U_{0})f_{1}(t+h^{*},U_{1},a)\mid T^{a}\geq t\right]<\infty$$
\end{condition}

\noindent The real value of the \mchr, that deviates from the \chr,  is derived in Theorem \ref{CH5th42}. 
\begin{theorem}\label{CH5th42}
	If the cause-effect relations of interest can be parameterized with SCM \eqref{CH5SCMsurv}, where 
	\begin{equation*}\lambda_{i}^{a}(t) = f_{0}(t,U_{0i})f_{1}(t,U_{1i},a),
	\end{equation*} and Condition \ref{CH5regcond} holds, then 
	\begin{equation*}
	\frac{\lim_{h\rightarrow 0}h^{-1}\mathbb{P}\left(T^{a} \in [t,t+h) \mid T^{a}\geq t \right)}{\lim_{h\rightarrow 0}h^{-1}\mathbb{P}\left(T^{0} \in [t,t+h) \mid T^{0}\geq t \right)} = \frac{\mathbb{E}\left[f_{0}(t,U_{0})f_{1}(t,U_{1},a)\mid T^{a}\geq t\right]}{\mathbb{E}\left[f_{0}(t,U_{0})\mid T^{0}\geq t\right]}=
	\end{equation*}
	\begin{equation*}
	\int f_{0}(t,U_{0})f_{1}(t,U_{1},a)\tfrac{\exp(-\Lambda^{a}(t,U_{0},U_{1})}{\int \exp(-\Lambda^{a}(t,U_{0},U_{1})dF_{(U_{0},U_{1})}}dF_{(U_{0},U_{1})}
	\left(\int f_{0}(t,U_{0})\tfrac{\exp(-\Lambda^{0}(t,U_{0}))}{\int \exp(-\Lambda^{0}(t,U_{0}))dF_{U_{0}}}dF_{U_{0}}\right)^{-1},
	\end{equation*} where $\Lambda^{a}(t,u_{0},u_{1})=\int_{0}^{t}f_{0}(s,u_{0})f_{1}(s,u_{1},a)ds$ and thus $\Lambda^{0}(t,u_{0})=\int_{0}^{t}f_{0}(s,u_{0})ds$.
\end{theorem}
\noindent The conditional expectations that determine the value of the \mchr equal weighted means of\\ $f_{0}(t,u_{0})f_{1}(t,u_{1},a)$ and $f_{0}(t,u_{0})$ with weights $\tfrac{\mathbb{P}(T^a\geq t\mid U_{0}=u_{0},U_{1}=u_{1})}{\mathbb{P}(T^{a}\geq t)}$
and $\tfrac{\mathbb{P}(T^0\geq t\mid U_{0}=u_{0})}{\mathbb{P}(T^{a}\geq t)}$ respectively. To develop our understanding of the bias when the \ohr (assuming no confounding) is used to estimate the \chr, we will first continue to study the bias as a result of frailty and heterogeneity separately in the next two subsections. 

\subsection{Causal effect homogeneity}\label{CH5frailnohet}

In the case of homogeneous multiplicative causal effects on the hazard, i.e.~$f_{1}(t, U_{1i}, a)$ $=$ $f_{1}(t,a)$, the ratio of the marginal hazard rates of individuals satisfying $T_{i}^{a}\geq t$ and of those $T_{i}^{0}\geq t$ equals $f_{1}(t,a)$ 
multiplied by a factor that depends on the difference in frailty distributions at time $t$ in those two populations as derived in Corollary \ref{CH5cor411}. 

\begin{corollary}\label{CH5cor411}
	If the cause-effect relations of interest can be parameterized with SCM \eqref{CH5SCMsurv}, where 
	\begin{equation*}\lambda_{i}^{a}(t) = f_{0}(t,U_{0i})f_{1}(t,a), 
	\end{equation*} and condition \ref{CH5regcond} holds 
	then 
	\begin{equation*}
	\frac{\lim_{h\rightarrow 0}h^{-1}\mathbb{P}\left(T^{a} \in [t,t+h) \mid T^{a}\geq t \right)}{\lim_{h\rightarrow 0}h^{-1}\mathbb{P}\left(T^{0} \in [t,t+h) \mid T^{0}\geq t \right)} = \frac{\mathbb{E}\left[f_{0}(t,U_{0})\mid T^{a}\geq t\right]}{\mathbb{E}\left[f_{0}(t,U_{0})\mid T^{0}\geq t\right]}f_{1}(t,a) = 
	\end{equation*}
	\begin{equation*}
\int f_{0}(t,U_{0})\tfrac{\exp(-\int_{0}^{t}f_{0}(s,U_{0})f_{1}(a,s)ds)} {\int \exp(-\int_{0}^{t}f_{0}(s,U_{0})f_{1}(a,s)ds)dF_{U_{0}}}dF_{U_{0}}
	\left(\int f_{0}(t,U_{0})\tfrac{\exp(-\int_{0}^{t}f_{0}(s,U_{0})ds)}{\int \exp(-\int_{0}^{t}f_{0}(s,U_{0}))dF_{U_{0}}}dF_{U_{0}}\right)^{-1}f_{1}(t,a).
	\end{equation*}
\end{corollary} \noindent The conditional expectation $\mathbb{E}\left[f_{0}(t,U_{0})\mid T^{a}\geq t\right]$ can be seen as a reweighted mean of $f_{0}(t,u_{0})$. These weights equal $\frac{\mathbb{P}(T^a\geq t\mid U_{0}=u_{0})}{\mathbb{P}(T^a\geq t)}$ and over time increase for favourable values of $U_{0}$. As a result 
	$\mathbb{E}\left[f_{0}(t,U_{0})\mid T^{a}\geq t\right]$ decreases over time. If $\Lambda^{a}(u_{0},t)>\Lambda^{0}(u_{0},t)$, e.g.~when $\forall t:f_{1}(t,a)>1$, then the weights $\frac{\mathbb{P}(T^a\geq t\mid U_{0}=u_{0})}{\mathbb{P}(T^a\geq t)}$ decrease more rapidly with $u$ than the weights $\frac{\mathbb{P}(T^0\geq t\mid U_{0}=u_{0})}{\mathbb{P}(T^0\geq t)}$, leading to 
\begin{equation}
\frac{\mathbb{E}\left[f_{0}(t,U_{0})\mid T^{a}\geq t\right]}{\mathbb{E}\left[f_{0}(t,U_{0})\mid T^{0}\geq t\right]}<1.
\end{equation} On the contrary, when $\Lambda^{a}(u_{0},t)<\Lambda^{0}(u_{0},t)$, 
then
\begin{equation}
\frac{\mathbb{E}\left[f_{0}(t,U_{0})\mid T^{a}\geq t\right]}{\mathbb{E}\left[f_{0}(t,U_{0})\mid T^{0}\geq t\right]}>1.
\end{equation}
So in the case of effect homogeneity, at time $t$, the \mchr is larger than the \chr when $\Lambda^{a}(u_{0},t)<\Lambda^{0}(u_{0},t)$ and smaller when $\Lambda^{a}(u_{0},t)>\Lambda^{0}(u_{0},t)$. An example of the latter was showcased by \citet{Stensrud2017}, where a model with $f_{1}(t, a)=1.81^a$, $f_{0}(u_{0},t)=u_{0}\lambda_{0}(t)$ and compound Poisson distributed frailty $U_{0}$ could well explain the decrease of the effect of hormone replacement therapy on coronary heart disease in postmenopausal women over time as observed from an RCT by the Woman Health Initiative. Furthermore, as emphasized by \citet{Hernan2010} based on the same case-study, even when $f_{1}(t,a)$ is constant over time, the (bias of the) \mchr  
hazard is time-varying, so that the estimates can depend on the follow-up time of a study. 

For frailty models as presented by \citet{Stensrud2017} and \citet{Aalen2015}, where
$
f_{0}(t,U_{0i})= U_{0i}\lambda_{0}(t)
$, it has already been shown by \citet{Balan2020} that $\mathbb{E}\left[U_{0}\mid T\geq t, A=a\right]$ can be expressed in terms of the Laplace transform of the frailty $U_{0}$. Reasoning along the same lines,  $\mathbb{E}\left[U_{0}\mid T^{a}\geq t\right]$ is expressed in terms of the Laplace transform of the $U_{0}$ in Lemma \ref{CH5lemma1}. 
\begin{lemma}\label{CH5lemma1}
	If the cause-effect relations of interest can be parameterized with SCM \eqref{CH5SCMsurv}, where 
	\begin{equation*}f_{\lambda}(t,U_{0i},U_{1i},a) = U_{0i}\lambda_{0}(t)f_{1}(t,a), 
	\end{equation*} then 
	\begin{equation}
	\mathbb{E}\left[U_{0} \mid T^{a}\geq t\right]   
	=-\frac{\mathcal{L}_{U_{0}}^{'}(\int_{0}^{t}\lambda_{0}(s)f_{1}(a,s)ds)}{\mathcal{L}_{U_{0}}(\int_{0}^{t}\lambda_{0}(s)f_{1}(a,s)ds)},
	\end{equation} where $\mathcal{L}_{U_{0}}(c)=\mathbb{E}\left[\exp(-cU_{0})\right]$ with derivative $\mathcal{L}_{U_{0}}^{'}(c)$.
\end{lemma} As in this previous work \citep[Figure 5]{Balan2020}, we present examples where $U_{0}$ follows a gamma, inverse Gaussian or compound Poisson distribution respectively. The parameterizations, corresponding Laplace transforms and expressions for $\mathbb{E}[U_{0}\mid T^{a}\geq t]$ can be found in Appendix \ref{CH5app:LT}. To illustrate the selection bias, we consider a binary exposure and let 
$$ \lambda_{i}^{a}(t) = U_{0i}\lambda_{0}(t)c^{a},$$ where $\lambda_{0}(t)=\tfrac{t^{2}}{20}$, $U_{0}\sim \Gamma(\theta_{0}^{-1},\theta_{0}) $, $U_{0}\sim \text{IG}(1,\theta_{0}^{-1})$ or $U_{0}\sim \text{CPoi}(3\theta_{0}^{-1},\tfrac{1}{2},\tfrac{2}{3}\theta_{0})$ respectively, so that $\mathbb{E}[U_{0}]=1$ and $\text{var}(U_{0})=\theta_{0}$. 
By applying Lemma 
\ref{CH5lemma1}, $\mathbb{E}[U_{0}\mid T^{1}\geq t]$ and $\mathbb{E}[U_{0}\mid T^{0}\geq t]$ can be derived. The \mchr then follows from Corollary \ref{CH5cor411} (conditional hazard is monotone increasing). 
The expressions for these quantities are presented in Table \ref{CH5tab:1}.
\begin{table}[h!]		\captionsetup{width=\textwidth}
		\centering
	\small
	\tabcolsep=0.10cm
	\caption{ Conditional expectations and resulting \mchr when the \chr equals $c$ for different frailty distributions such that $\mathbb{E}[U_{0}]=1$ and $\text{var}(U_{0})=\theta_{0}$ (in absence of effect modification). \label{CH5tab:1}}
	\begin{tabular}{c||c|c|c} 
		& \multicolumn{3}{c}{\Gape[0.25cm][0.25cm]{$U_{0}\sim$}} \\ \cline{2-4}
		& \Gape[0.25cm][0.25cm]{$\Gamma(\theta_{0}^{-1},\theta_{0})$ }   & $\text{IG}(1,\theta_{0}^{-1})$                               & $\text{CPoi}\left(3\theta_{0}^{-1},\tfrac{1}{2},\tfrac{2}{3}\theta_{0}\right)$ \\ \hline
		$\mathbb{E}[U_{0}\mid T^{1}\geq t]$                                                                                                                                                & \Gape[0.25cm][0.25cm]{$\frac{60}{60+ct^{3}\theta_{0}}$}        & $\sqrt{\frac{30}{30+ct^{3}\theta_{0}}}$                      & $90^{\tfrac{3}{2}}\left(90+ct^3\theta_{0}\right)^{-\tfrac{3}{2}}$   \\ \hline
		$\mathbb{E}[U_{0}\mid T^{0}\geq t]$                                                                                                                                                & \Gape[0.25cm][0.25cm]{$\frac{60}{60+t^{3}\theta_{0}}$}         &                   $\sqrt{\frac{30}{30+t^{3}\theta_{0}}}$ &          $90^{\tfrac{3}{2}}\left(90+t^3\theta_{0}\right)^{-\tfrac{3}{2}}$                                                                        \\ \hline
		\Gape[0.25cm][0.25cm]{$\tfrac{\lim_{h\rightarrow 0}h^{-1}\mathbb{P}\left(T^{1} \in (t,t+h) \mid T^{1}\geq t \right)}{\lim_{h\rightarrow 0}h^{-1}\mathbb{P}\left(T^{0} \in [t,t+h) \mid T^{0}\geq t \right)}$} & $1+\frac{60(c-1)}{60+ct^{3}\theta_{0}}$ & $c \frac{\sqrt{30+t^3\theta_{0}}}{\sqrt{30+ct^3\theta_{0}}}$ & $c \left(\frac{90+t^3\theta_{0}}{90+ct^3\theta_{0}}\right)^{\tfrac{3}{2}}$          
	\end{tabular}
\end{table} 

\begin{figure}[h!]
	\centering
	\captionsetup{width=\textwidth}
	\includegraphics[width=\textwidth]{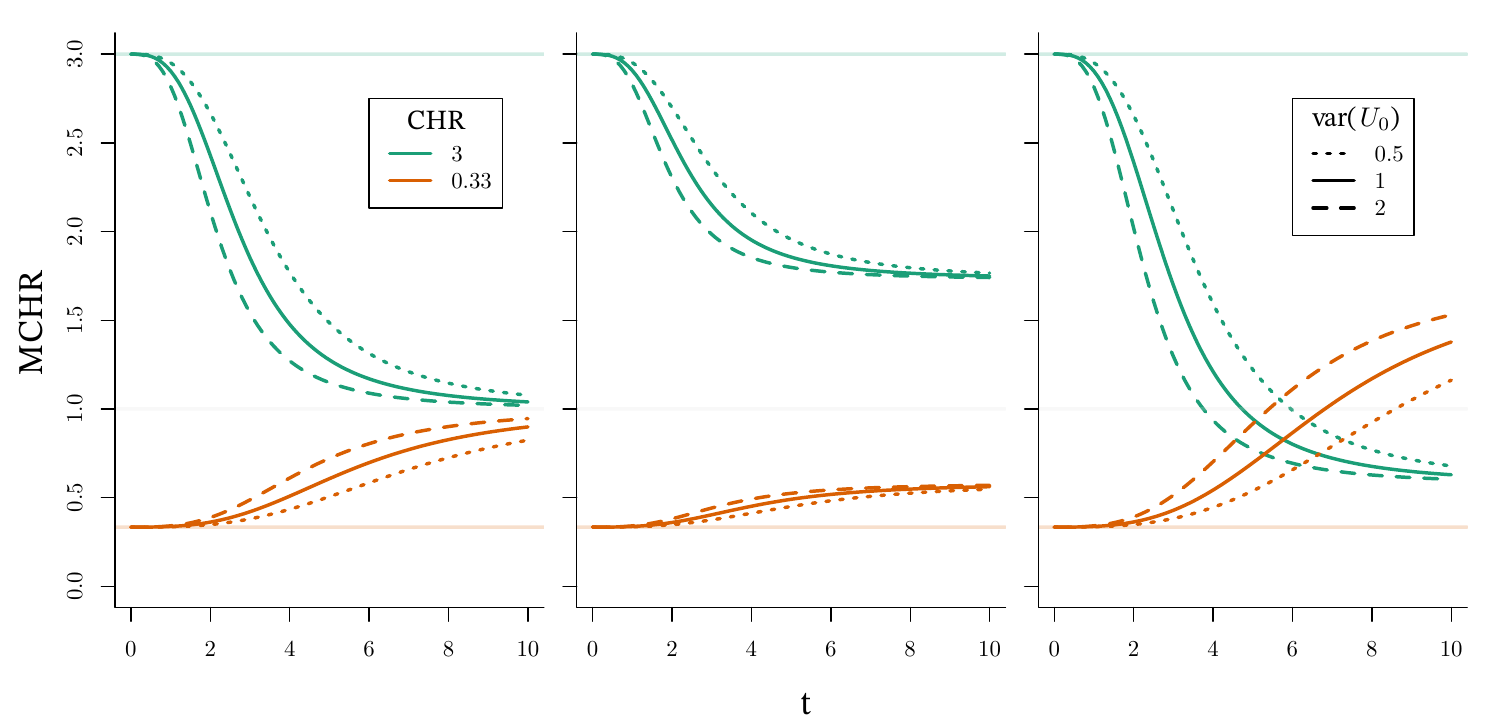} 
	\caption{\mchr  over time when $\lambda^{a}_{i}(t)=U_{0i} \tfrac{t^{2}}{20}c^{a}$ for $c=3$ (green) and $c=\tfrac{1}{3}$ (orange), when $U_{0}$ follows a gamma (left), inverse Gaussian (middle) or compound Poisson (right) distribution with variance $0.5$ (dotted), $1$ (solid) or $2$ (dashed) respectively.  }\label{CH5Frail}	
\end{figure}

How the \mchr
deviates from $c$ over time for $c\in\{\tfrac{1}{3}, 3\}$,  and $\theta \in \{0.5, 1, 2\}$  is visualized in Figure \ref{CH5Frail}. 
For both $c=3$ and $c=\tfrac{1}{3}$ the selection of individuals that survive time $t$ result in a \mchr that evolves in the opposite direction of the causal effect, towards $1$, $\sqrt{c}$ and $\sqrt{c^{-1}}$ respectively. For the case of a compound Poisson frailty, the logarithm of this latter limit is even opposite to the sign of the logarithm of
the \chr due to the nonsusceptible individuals. For all types of frailty, the higher the variance of $U_{0}$, the larger the difference between the \mchr and the \chr. For comparison we have also presented the survival curves of $T^{1}$ and $T^{0}$ in Figure \ref{CH5Surv1} in Appendix \ref{CH5app:surv} for the setting where $\theta_{0}=1$. Note that for an RCT $T^{a}\overset{d}{=} T \mid A=a$ by \eqref{CH5rand} and causal consistency.

\subsection{Causal effect heterogeneity in the absence of frailty}\label{CH5hetnofrail}

Before we return to the general case presented in Theorem \ref{CH5th42}, let's consider the presence of effect heterogeneity in the absence of frailty, i.e.~
\begin{equation*}
f_{\lambda}(t,U_{0i},U_{1i},a) = \lambda_{0}(t)f_{1}(t,U_{1i},a).
\end{equation*} If the \chr,  $\mathbb{E}[f_{1}(t,U_{1i},a)]$, is equal for all $t$, the \mchr is not as over time among the exposed individuals those that `benefit' more are more likely to survive. The effect of this selection on the \mchr over time is formalized in Corollary \ref{CH5cor421}.  
\begin{corollary}\label{CH5cor421}
	If the cause-effect relations of interest can be parameterized with SCM \eqref{CH5SCMsurv}, where 
	\begin{equation*}\lambda_{i}^{a}(t) = \lambda_{0}(t)f_{1}(t,U_{1i},a), 
	\end{equation*} and Condition \ref{CH5regcond} holds then 
	\begin{align*}
	\frac{\lim_{h\rightarrow 0}h^{-1}\mathbb{P}\left(T^{a} \in [t,t+h) \mid T^{a}\geq t \right)}{\lim_{h\rightarrow 0}h^{-1}\mathbb{P}\left(T^{0} \in [t,t+h) \mid T^{0}\geq t \right)} &= \mathbb{E}\left[f_{1}(t,U_{1},a)\mid T^{a}\geq t\right]\\	&= \int f_{1}(t,U_{1},a)\tfrac{\exp(-\Lambda^{a}(t,U_{1}))}{\int \exp(-\Lambda^{a}(t,U_{1}))dF_{U_{1}}} dF_{U_{1}},
	\end{align*} where $\Lambda^{a}(t,u_{1})=\int_{0}^{t}\lambda_{0}(s)f_{1}(s,u_{1},a)ds$.
\end{corollary} \noindent The \mchr thus equals $\mathbb{E}\left[f_{1}(t,U_{1},a)\mid T^{a}\geq t\right]$, which is smaller than $\mathbb{E}\left[f_{1}(t,U_{1},a)\right]$
as more weight is placed on lower values of $f_{1}(t,u_{1},a)$ that correspond to higher $\mathbb{P}(T^a\geq t\mid U_{1}=u_{1})$. Besides the selection of frailty factors, the selection of individual modifiers can also lead to selection bias of the estimated hazard ratio. In contrast to $\mathbb{E}[f_{0}(t,U_{0})]$ in the previous section, 
$\mathbb{E}\left[f_{1}(t,U_{1},a)\mid T^{a}\geq t\right]$ decreases over time irrespective of whether the exposure is beneficial or harming on average. For a hypothetical setting without frailty but with effect heterogeneity, the \chr at $t$ is thus systematically underestimated when using the \mchr, so the exposure seems more `beneficial' than it is. Such a difference has only been explained by the presence of frailty \citep{Hernan2010, Stensrud2017}. For binary exposures and SCMs that further restrict  $
f_{1}(t,u_{1},a)= u_{1}^a$, 
by Lemma \ref{CH5lemma1} (as $\lambda_{i}^{1}=U_{1i}\lambda_{0}(t)$), 
\begin{equation*}
\mathbb{E}\left[U_{1}\mid T^{1}\geq t\right]=-\frac{\mathcal{L}_{U_{1}}^{'}(\int_{0}^{t}\lambda_{0}(s)ds)}{\mathcal{L}_{U_{1}}(\int_{0}^{t}\lambda_{0}(s)ds)}. 
\end{equation*} Let $\lambda_{0}=\tfrac{t^{2}}{20}$, $U_{1}\sim \Gamma(\tfrac{c}{\theta_{1}},\tfrac{\theta_{1}}{c}) $, $U_{1}\sim \text{IG}(c,\tfrac{c^3}{\theta_{1}})$ or $U_{1}\sim \text{CPoi}(3\tfrac{c^2}{\theta_{1}},\tfrac{1}{2},\tfrac{2\theta_{1}}{3c})$  respectively such that $\mathbb{E}[U_{1}]=c$ and $\text{var}(U_{1})=\theta_{1}$. By applying  Lemma \ref{CH5lemma1} for $a=1$, we can derive $\mathbb{E}[U_{1}\mid T^{1}\geq t]$, which by Corollary \ref{CH5cor421} (conditional hazard is monotone increasing) equal the \mchrs, and are presented in Table \ref{CH5tab:2}. Additionally, we derived $\mathbb{E}[U_{1}\mid T^{1}\geq t]$ for a setting where the multiplicative hazard effect modifier $U_{1}$ equals $\mu_{1}$ ($<1$, for individuals that benefit) with probability $p_{1}$, $\mu_{2}$ ($>1$, for individuals that are harmed) with probability $p_{2}$ or $1$ (for individuals that are not affected). We define this distribution as the Benefit-Harm-Neutral, $\text{BHN}(p_{1},\mu_{1},p_{2},\mu_{2})$, distribution. 
\begin{table}[h!]
	\centering
	\small
	\tabcolsep=0.10cm
	\caption{ Conditional expectation of the individual effect modifier $U_{1}$ when the \chr equals $c$ for different modifier distributions such that $\mathbb{E}[U_{1}]=c$ and $\text{var}(U_{0})=\theta_{1}$ (in absence of frailty).\label{CH5tab:2}}
		\begin{tabular}{c||c|c|c|c}
			& \multicolumn{4}{c}{\Gape[0.25cm][0.25cm]{$U_{1}\sim$}} \\ \cline{2-5}
			& \Gape[0.25cm][0.25cm]{$\Gamma(\tfrac{c}{\theta_{1}},\tfrac{\theta_{1}}{c})$} & $\text{IG}\left(c,\tfrac{c^3}{\theta_{1}}\right)$              & $\text{CPoi}\left(3\tfrac{c^2}{\theta_{1}},\tfrac{1}{2},\tfrac{2\theta_{1}}{3c}\right)$ & $\text{BHN}(p_{1},\mu_{1},p_{2},\mu_{2})$\\ \hline
			$\mathbb{E}[U_{1}\mid T^{1}\geq t]$ & $\frac{60 c^{2}}{\theta_{1}t^{3}+60c}$ 
			& $\frac{c\sqrt{30c}}{\sqrt{ t^{3}\theta_{1}+30c}}$ & $c\left(\frac{\theta_{1}t^{3}}{90c}+1\right)^{-\tfrac{3}{2}}$     &
			\Gape[0.25cm][0.25cm]{$\frac{p_{1}\mu_{1}+p_{2}\mu_{2}e^{\tfrac{-t^3(\mu_{2}-\mu_{1})}{60}}+(1-p_{1}-p_{2})e^{\tfrac{-t^3(1-\mu_{1})}{60}}}{p_{1}+p_{2}e^{\tfrac{-t^3(\mu_{2}-\mu_{1})}{60}}+(1-p_{1}-p_{2})e^{\tfrac{-t^3(1-\mu_{1})}{60}}}$}
	\end{tabular}
\end{table}

\begin{figure}[h!]
	\centering
	\captionsetup{width=\textwidth}
	\includegraphics[width=\textwidth]{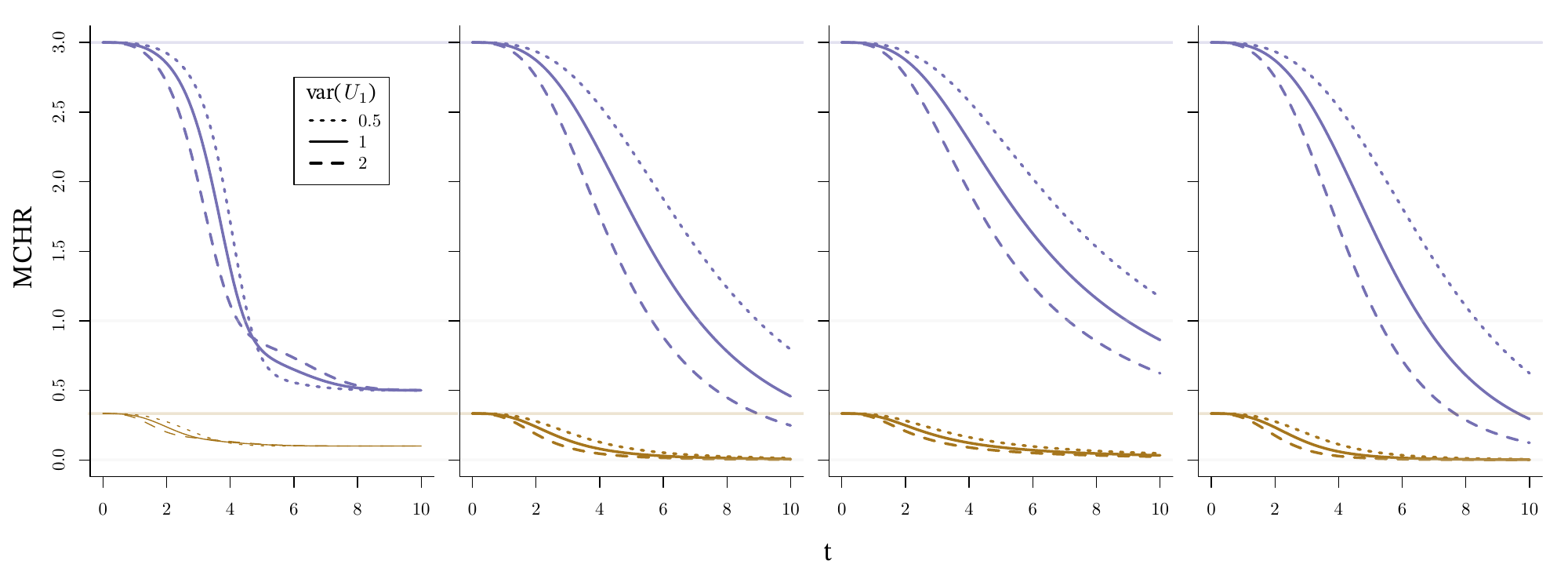} 
	\caption{\mchr  over time when $\lambda^{a}_{i}(t)=\tfrac{t^{2}}{20}(U_{1i})^{a}$, when $U_{1}$ follows a BHN, gamma, inverse Gaussian or compound Poisson distribution (from left to right) with expectation $3$ (blue) or $\tfrac{1}{3}$ (brown) and variance $0.5$ (dotted), $1$ (solid) or $2$ (dashed) respectively.  }\label{CH5Het}	
\end{figure}

For $\mathbb{E}[U_{1}] \in \left \{\tfrac{1}{3}, 3 \right \}$, and $\theta_{1} \in \{0.5, 1, 2\}$ the evolution of the conditional expectation is shown in Figure \ref{CH5Het} for all four effect-modifier distributions. For the $\text{BHN}$ distribution, when $\mathbb{E}[U_{1}]$ equals $3$ and $\tfrac{1}{3}$, we fix $p_{1}=0.05$, $\mu_{1}=0.5$ and $p_{1}=0.9$, $\mu_{1}=0.1$ respectively. Expressions for $p_{2}$ and $\mu_{2}$ such that $\mathbb{E}[U_{1}]=\mu$ and $\text{var}(U_{1})=\theta_{1}$ can be found in Appendix \ref{CH5discrete}. 
 When the exposure is in expectation harming $(\mathbb{E}[U_{1}]=3)$, for all settings considered there is a point in time that the \mchr drops below $1$. For the continuous distributions, the \mchr won't stop decreasing. The decreases for the gamma and compound Poisson settings are very similar, while for the inverse Gaussian setting, this goes a bit slower. For the discrete setting, the marginal causal hazard converges to $\mu_{1}$ of $0.5$ and $0.1$, respectively. Again, as in the previous subsection, the higher the variability of the latent variable, the faster the \mchr deviates from the \chr. Only for the discrete effect modifier, the lines cross for the different variances for  $\mathbb{E}[U_{1}=3]$, but this is the result of different fractions of individuals that are not affected by the exposure (as the mean and variance are coupled). 

\subsection{Causal effect heterogeneity in the presence of frailty}

In the general case where effect heterogeneity and frailty are present, both principles affect the value of the \mchr. By Theorem \ref{CH5th42}, the ratio evolves as\\  $\frac{\mathbb{E}\left[f_{0}(t,U_{0})f_{1}(t,U_{1},a)\mid T^{a}\geq t\right]}{\mathbb{E}\left[f_{0}(t,U_{0})\mid T^{0}\geq t\right]}$. The numerator depends on the joint distribution of $U_{0}$ and $U_{1}$. For illustration, we again consider a binary exposure and let
\begin{equation*}f_{\lambda}(t,U_{0i},U_{1i},a) = U_{0i}(U_{1i})^{a}\lambda_{0}(t)f_{1}(t,a), 
\end{equation*} such that the \mchr equals $\frac{\mathbb{E}\left[U_{0}U_{1}\mid T^{1}\geq t\right]}{\mathbb{E}\left[U_{0}\mid T^{0}\geq t\right]}f_{1}(1,t)$ and, by Lemma \ref{CH5lemma1}, can be derived from the Laplace transforms of $U_{0}U_{1}$ and $U_{0}$ respectively.

\subsubsection{Independent $U_{0}$ and $U_{1}$}
In the case of independence, the Laplace transform of the product equals $\mathbb{E}[\mathcal{L}_{U_{0}}(c U_{1})]$, which generally does not adopt a tractable form. The case with a discrete effect modifier, introduced in Section \ref{CH5hetnofrail}, forms an exception. If $U_{1}\sim \text{BHN}(p_{1},\mu_{1},p_{2},\mu_{2}$), then \begin{equation}
\mathcal{L}_{U_{0}U_{1}}(c)= p_{1}\mathcal{L}_{U_{0}}(\mu_{1}c)+p_{2}\mathcal{L}_{U_{0}}(\mu_{2}c)+(1-p_{1}-p_{2})\mathcal{L}_{U_{0}}(c).
\end{equation} For the running example where $f_{1}(t,a)=1$, $\lambda_{0}=\tfrac{t^{2}}{20}$, $U_{0}\sim \Gamma(\theta_{0}^{-1},\theta_{0}) $, $U_{0}\sim \text{IG}(1,\theta_{0}^{-1})$ or $U_{0}\sim \text{CPois}(3\theta_{0}^{-1},\tfrac{1}{2},\tfrac{2}{3}\theta_{0})$ for $\theta_{0}\in\{0.5, 1, 2\}$, the expressions for $\mathbb{E}[U_{0}U_{1}\mid T^{1}\geq t]$ are presented in Table \ref{CH5tab:3} (where $p_{3}=1-p_{1}+p_{2}$ and $\mu_{3}=1$). As the $\mathbb{E}[U_{0}\mid T^{0}\geq t]$ are independent of the $U_{1}$ distribution these expectations are the same as presented in Table \ref{CH5tab:1}. 

\begin{figure}[H]
	\centering
	\captionsetup{width=\textwidth}
	\includegraphics[width=0.9\textwidth]{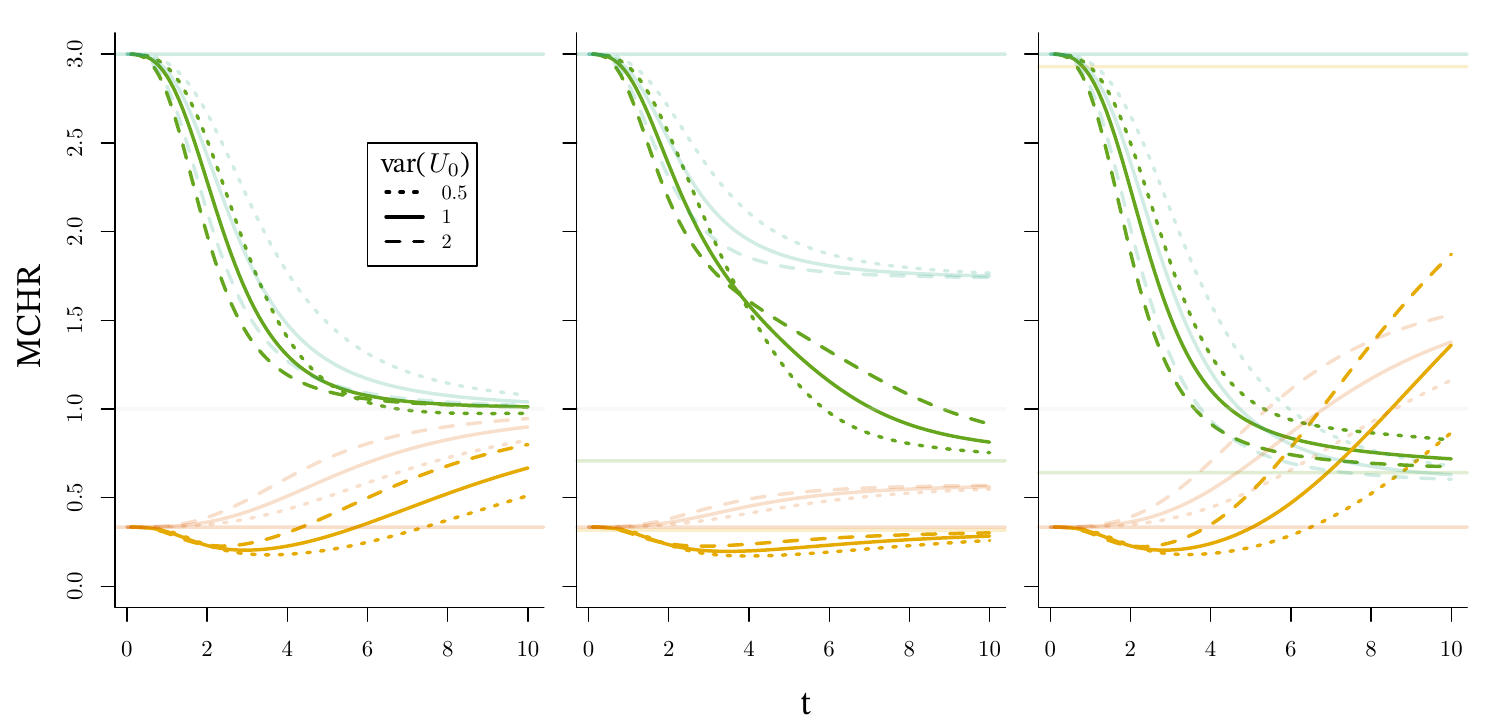} 
	\caption{\mchr  over time when $\lambda^{a}_{i}(t)=U_{0i}(U_{1i})^{a} \tfrac{t^{2}}{20}$ for a unit-variance BHN distributed $U_{1}$ with $\mathbb{E}[U_{1}]=3$ (opaque green) and $\mathbb{E}[U_{1}]=\tfrac{1}{3}$ (opaque orange), when $U_{0}$ follows a gamma (left), inverse Gaussian (middle) or compound Poisson (right) distribution with variance $0.5$ (dotted), $1$ (solid) or $2$ (dashed) respectively. For comparison, the lines presented in Figure \ref{CH5Frail} 
		are represented by transparent lines. }\label{CH5FrailDisc}	
\end{figure}

\begin{figure}[h!]
	\centering
	\captionsetup{width=\textwidth}
	\includegraphics[width=0.9\textwidth]{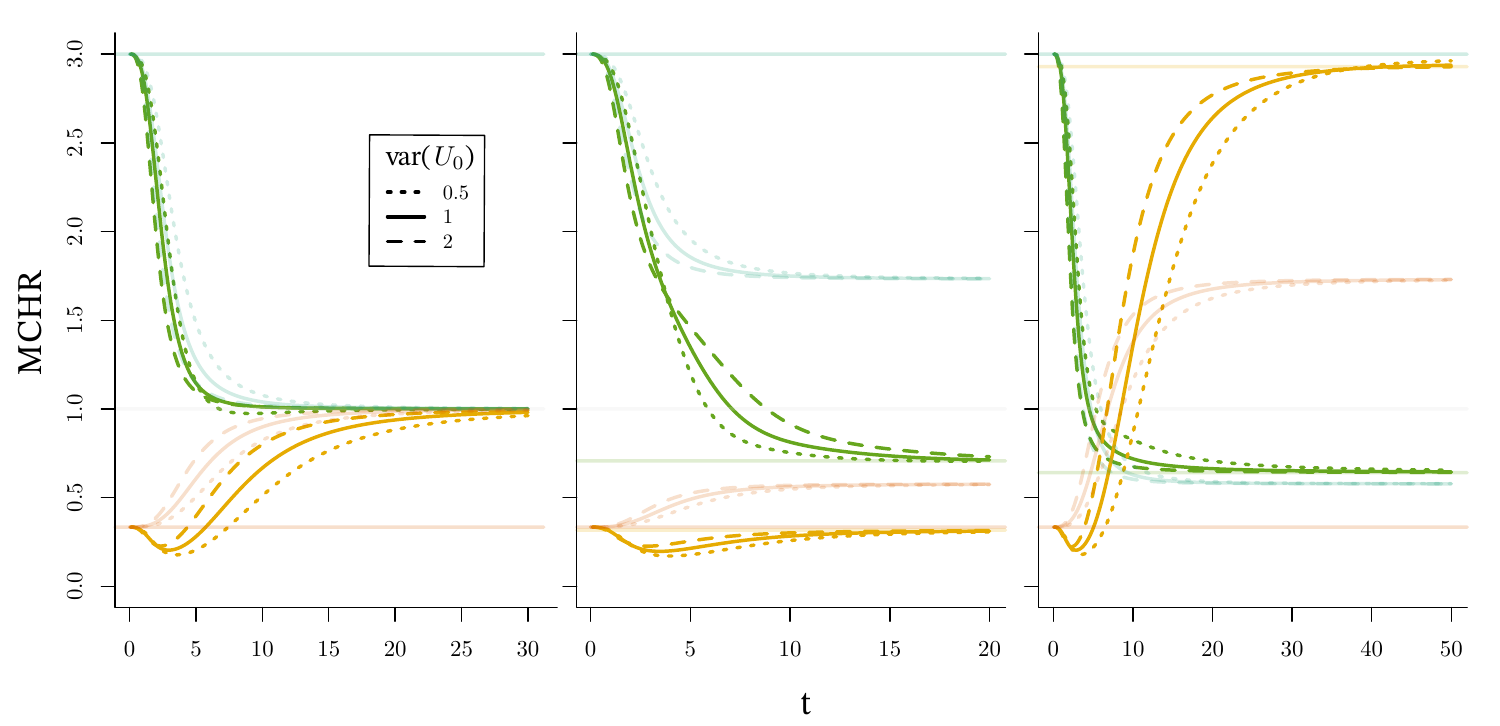}
	\caption{\mchr  over time when $\lambda^{a}_{i}(t)=U_{0i}(U_{1i})^{a} \tfrac{t^{2}}{20}$ for a unit-variance BHN distributed $U_{1}$ with $\mathbb{E}[U_{1}]=3$ (opaque green) and $\mathbb{E}[U_{1}]=\tfrac{1}{3}$ (opaque orange), when $U_{0}$ follows a gamma (left), inverse Gaussian (middle) or compound Poisson (right) distribution with variance $0.5$ (dotted), $1$ (solid) or $2$ (dashed) respectively. For comparison, the lines presented in Figure \ref{CH5Frail} 
		as well as the limits of the \mchrs are represented by transparent lines. }\label{CH5FrailDiscL}	
\end{figure}

\begin{table}[h!]
	\centering
\small
\tabcolsep=0.10cm
	\caption{ Conditional expectation $\mathbb{E}[U_{0}U_{1}\mid T^{1}\geq t]$ and its limiting forms when $U_{1}$ follows a $\text{BHN}(p_{1},\mu_{1},p_{2},\mu_{2})$ distribution) for different unit-expectations frailty distributions with $\text{var}(U_{0})=\theta_{0}$ while $U_{1} \independent U_{0}$. For comparison $\mathbb{E}[U_{0} \mid T^{0}\geq t]$, presented before in Table \ref{CH5tab:1}, is also presented  \label{CH5tab:3}}. 
	\begin{tabular}{c||c|c} 
		\Gape[0.25cm][0.25cm]{$U_{0}\sim$}                                                         & $\mathbb{E}[U_{1}U_{0}\mid T^{1}\geq t]$                                                                                                                    & $\mathbb{E}[U_{0}\mid T^{0}\geq t]$                                  \\ \hline
		$\Gamma(\theta_{0}^{-1},\theta_{0})$                                & \Gape[0.25cm][0cm]{$\frac{\sum_{i=1}^{3} p_{i}\mu_{i}\left(1+\tfrac{\theta_{0}t^3}{60}\mu_{i}\right)^{-(1+\theta_{0}^{-1})}}{\sum_{i=1}^{3} p_{i}(1+\tfrac{\theta_{0}t^3}{60}\mu_{i})^{-\theta_{0}^{-1}}}$} & $\frac{60}{60+t^{3}\theta_{0}}$                                  \\ 
		& \Gape[0cm][0.25cm]{= $(\theta_{0}\tfrac{t^3}{60})^{-1} + o\left(t^{-3}\right)$}& \\  \hline
		$\text{IG}(1,\theta_{0}^{-1})$                                      &      \Gape[0.25cm][0cm]{ $\frac{\sum_{i=1}^{3} p_{i}\mu_{i}\left(2\theta_{0}\tfrac{t^3}{60}\mu_{i}\right)^{-\tfrac{1}{2}}\exp \left(\theta_{0}^{-1}\left(1-\sqrt{1+\theta_{0}2\tfrac{t^{3}}{60}\mu_{i} }\right)\right)}{\sum_{i=1}^{3} p_{i}\exp \left(\theta_{0}^{-1}\left(1-\sqrt{1+\theta_{0}2\tfrac{t^{3}}{60}\mu_{i} }\right)\right)}$ } &  $\sqrt{\frac{30}{30+t^{3}\theta_{0}}}$                       \\
		& \Gape[0cm][0.25cm]{$= \sqrt{\mu_{1}}\sqrt{\tfrac{30}{t^{3}\theta_{0}}} + o\left(t^{-\tfrac{3}{2}}\right)$} &   \\\hline 
		$\text{CPoi}\left(3\theta_{0}^{-1},\tfrac{1}{2},\tfrac{2}{3}\theta_{0}\right)$ &                                                            \Gape[0.25cm][0cm]{$\frac{\sum_{i=1}^{3} p_{i}\mu_{i} \left(\frac{3}{3+2\theta_{0}\mu_{i}\tfrac{t^3}{60}}\right)^{\tfrac{3}{2}}\exp\left(3\theta_{0}^{-1}\sqrt{\frac{3}{3+2\theta_{0}\mu_{i}\tfrac{t^3}{60}}}-1\right)}{\sum_{i=1}^{3} p_{i}\exp\left(3\theta_{0}^{-1}\sqrt{\frac{3}{3+2\theta_{0}\mu_{i}\tfrac{t^3}{60}}}-1\right)}$} & $90^{\tfrac{3}{2}}\left(90+t^3\theta_{0}\right)^{-\tfrac{3}{2}}$\\
		& = \Gape[0cm][0.25cm]{$90^{\tfrac{3}{2}}\left(\theta_{0}t^{3}\right)^{-\tfrac{3}{2}}\sum_{i=1}^{3} \frac{p_{i}}{\sqrt{\mu_{i}}} + o\left(t^{-\tfrac{9}{2}}\right)$} & 

	\end{tabular}%
\end{table} The \mchr and its limit can be derived from the expressions in Table \ref{CH5tab:3} by applying Theorem \ref{CH5th42} (conditional hazard is monotone increasing). Interestingly, for gamma frailty, this limit remains $1$. For the inverse Gaussian frailty, the effect heterogeneity drastically change the limit from $\sqrt{\mathbb{E}[U_{1}]}$ to $\sqrt{\mu_{1}}$, which is always less than 1. Finally, for compound Poisson frailty the limit changes from $\sqrt{\mathbb{E}[U_{1}]}^{-1}$ to $\frac{p_{1}}{\sqrt{\mu_{1}}}+\frac{p_{2}}{\sqrt{\mu_{2}}} + (1-p_{1}-p_{2})$. The evolution of $\frac{\mathbb{E} \left[ U_{0}U_{1}\mid T^{1}\geq t \right] }{\mathbb{E} \left[ U_{0}\mid T^{0}\geq t \right] }$ over time is visualized in Figure \ref{CH5FrailDisc} for $U_{1}\sim \text{BHN}(0.05,0.5,0.82,3.5)$, such that $\mathbb{E}[U_{1}]=3$ and $\text{var}(U_{1})=1$, and for $U_{1}\sim \text{BHN}(0.9,0.1,0.03,6.0)$, such that $\mathbb{E}[U_{1}]=\tfrac{1}{3}$ and $\text{var}(U_{1})=1$.  In the case the \chr is larger than $1$, the selection of less susceptible individuals (frailty) that are harmed less (effect modifier) in the exposed world, both cause the \mchr to be smaller than the \chr. The \mchr decreases faster in the presence of effect heterogeneity. This explains the observation by \citet{Stensrud2017}, \textit{``Interestingly, the magnitude of frailty bias is larger when a heterogeneous treatment effect is included"}, for a simulation with frailty and random individual hazard ratios such that  $\mathbb{E}[\lambda_{i}^{1}(t)]=1.81>1$.
For the gamma and compound Poisson frailty examples this effect is relatively small as $\mathbb{E}[U_{1}U_{0}\mid T^{1}\geq t]$ is quite similar to $\mathbb{E}[U_{0}\mid T^{1}\geq t]$ as presented in Table \ref{CH5tab:1} (for the presented $p_{1}, \mu_{1}, p_{2}$ and $\mu_{2}$). However, for inverse Gaussian frailty, the \mchr deviates much more from the \chr in the presence of effect heterogeneity. In Figure \ref{CH5FrailDiscL}, the evolution of the \mchr is presented for a longer timescale, and the limits become apparent.  If the \chr is smaller than $1$, then the selection of less susceptible individuals (frailty) in the unexposed world and the selection of individuals that benefit more (effect modifier) in the exposed world have opposite effects on the \mchr. For this case of discrete effect modifiers, the \mchr first decreases by selecting individuals with more beneficial modifiers and later increases (above the \chr) when the frailty selection effect prevails. For the examples presented, the fraction $p_{1}=0.9$ of the population with $\mu_{1}=0.1$ are expected to survive so that over time the \mchr will resemble the \mchr in the absence of effect heterogeneity for this subpopulation (with the \chr equal to $0.1$). The limit for gamma frailty is still one, so the \mchr deviates less from the \chr due to the two opposed selection effects. The bias is strongly reduced for the inverse Gaussian frailty as the limit $\sqrt{0.1}$ is close to the actual \chr. Finally, for the compound Poisson frailty, the \mchr with effect heterogeneity crosses the \mchr in the absence of effect heterogeneity as the frailty bias is larger for a \chr of $0.1$ compared to one of $\tfrac{1}{3}$. In summary, the bias for the \chr can further increase in the presence of effect heterogeneity, stressing the issues regarding the causal interpretation of \ohrs (assuming no confounding). However, for beneficial exposures, the frailty bias can reduce in the presence of effect heterogeneity (e.g., Inverse Gaussian frailty), illustrating that there might be settings where the \mchr is close to the \chr.

\subsubsection{Dependent $U_{0}$ and $U_{1}$}
\begin{figure}[h!]
	\centering
	\captionsetup{width=\textwidth}
	\includegraphics[width=0.9\textwidth]{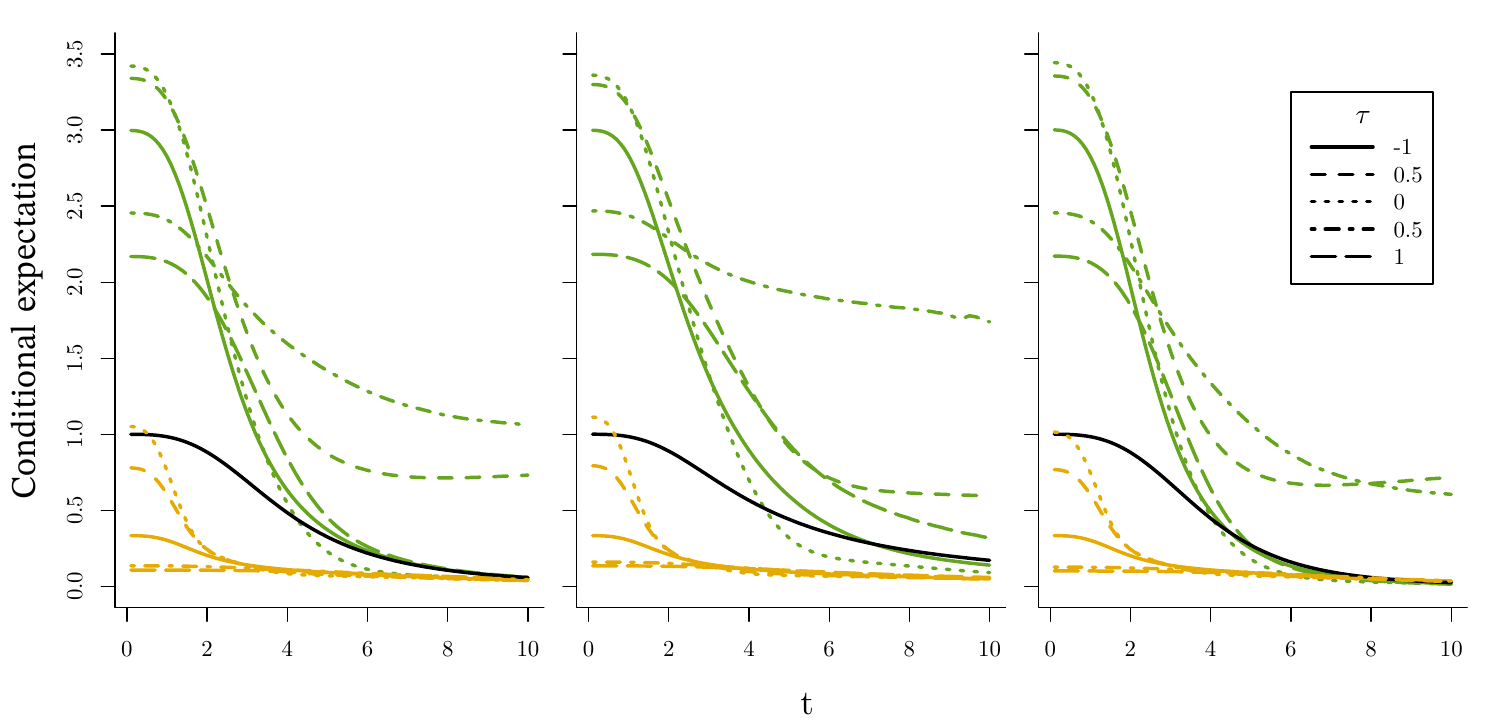} 			
	\caption{$\mathbb{E}[U_{0}U_{1}\mid T^{1}\geq t]$ over time for $\lambda^{a}_{i}(t)=U_{0i}(U_{1i})^{a} \tfrac{t^{2}}{20}$, for a unit-variance BHN distributed $U_{1}$ with $\mathbb{E}[U_{1}]=3$ (green) or $\mathbb{E}[U_{1}]=\tfrac{1}{3}$ (orange), $U_{0}$ follows a gamma (left), inverse Gaussian (middle) or compound Poisson (right) distribution and where the joint distribution of $U_{0}$ and $U_{1}$ follows from a Gaussian copula with varying Kendall's $\tau$ correlation coefficient (see legend). Furthermore, the evolution of $\mathbb{E}[U_{0}\mid T^{0}\geq t]$ is presented (black). }\label{CH5Fig5a}	
\end{figure}

\begin{figure}[h!]
	\centering
	\captionsetup{width=\textwidth}
	\includegraphics[width=0.9\textwidth]{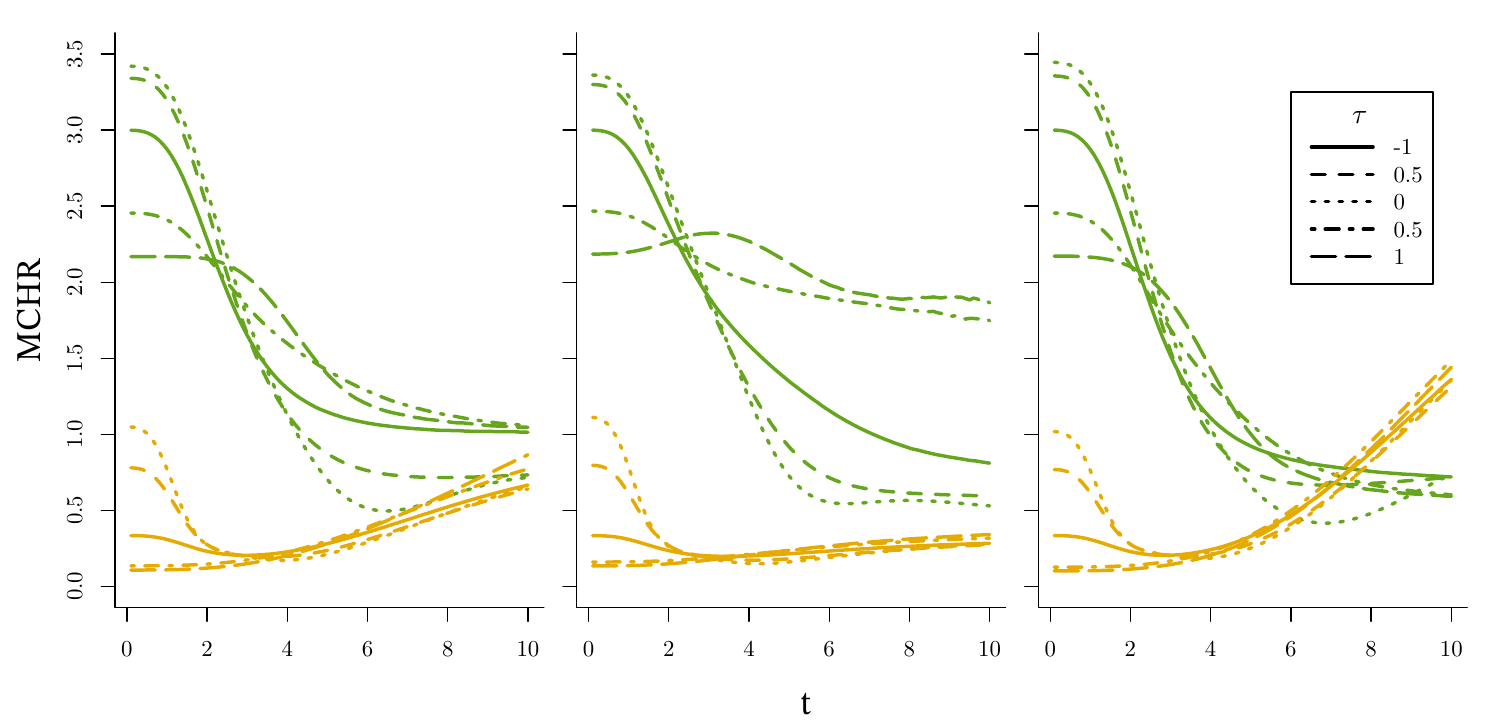}
	
	\caption{\mchr  over time for $\lambda^{a}_{i}(t)=U_{0i}(U_{1i})^{a} \tfrac{t^{2}}{20}$, for a unit-variance BHN distributed $U_{1}$ with $\mathbb{E}[U_{1}]=3$ (green) or $\mathbb{E}[U_{1}]=\tfrac{1}{3}$ (orange), $U_{0}$ follows a gamma (left), inverse Gaussian (middle) or compound Poisson (right) distribution and where the joint distribution of $U_{0}$ and $U_{1}$ follows from a Gaussian copula with varying Kendall's $\tau$ correlation coefficient (see legend). }\label{CH5Fig5}	
\end{figure}

In case the multiplicative effect of the exposure on the hazard of  susceptible individuals is expected to be higher or lower than for less susceptible individuals, the distribution of $U_{0}U_{1}$ will be less or more variable than when the latent variables are independent. Every bivariate joint distribution function can be written using the marginal distribution functions and a copula $C$ \citep{Sklar1959}. As such, 
\begin{equation*}
F_{(U_{0},U_{1})}(u_{0},u_{1})=C\left(F_{U_{0}}(u_{0}), F_{U_{1}}(u_{1})\right)  
\end{equation*} and the Kendall's $\tau$ correlation coefficient of $U_{0}$ and $U_{1}$ can be written as a function of the Copula $C$ \citep{Nelsen2006}. To study how the dependence can affect the \mchr for the setting presented in Figure \ref{CH5FrailDisc}, we use a Gaussian copula \begin{equation*}
C(x,y)=\Phi_{2,\rho}(\Phi^{-1}(x), \Phi^{-1}(y)),
\end{equation*} where $\Phi$ and $\Phi_{2,\rho}$ are the standard normal and bivariate normal with correlation $\rho$, cumulative distribution functions, respectively. For $\rho \in \{-1, \sin(-\tfrac{\pi}{4}), 0, \sin(\tfrac{\pi}{4}), 1\}$ (such that $\tau \in \{-1, -0.5, 0, 0.5, 1\}$) and $\text{var}(U_{0})=1$, $\mathbb{E}\left[U_{0}U_{1}\mid T^{1}\geq t\right]$ is derived empirically from simulations and are presented in Figure \ref{CH5Fig5a}.   All programming codes used for this work can be found online at \url{https://github.com/RAJP93/CHR}. The results were very similar when using a Frank, Clayton or Gumbel copula instead of the Gaussian copula. For reference, the evolution of  $\mathbb{E}\left[U_{0} \mid T^{0}\geq t\right]$ over time is also presented in Figure \ref{CH5Fig5a}. The corresponding \mchrs are presented in Figure \ref{CH5Fig5}.  

\begin{figure}[h!]
	\centering
	\captionsetup{width=\textwidth}
	\includegraphics[width=0.9\textwidth]{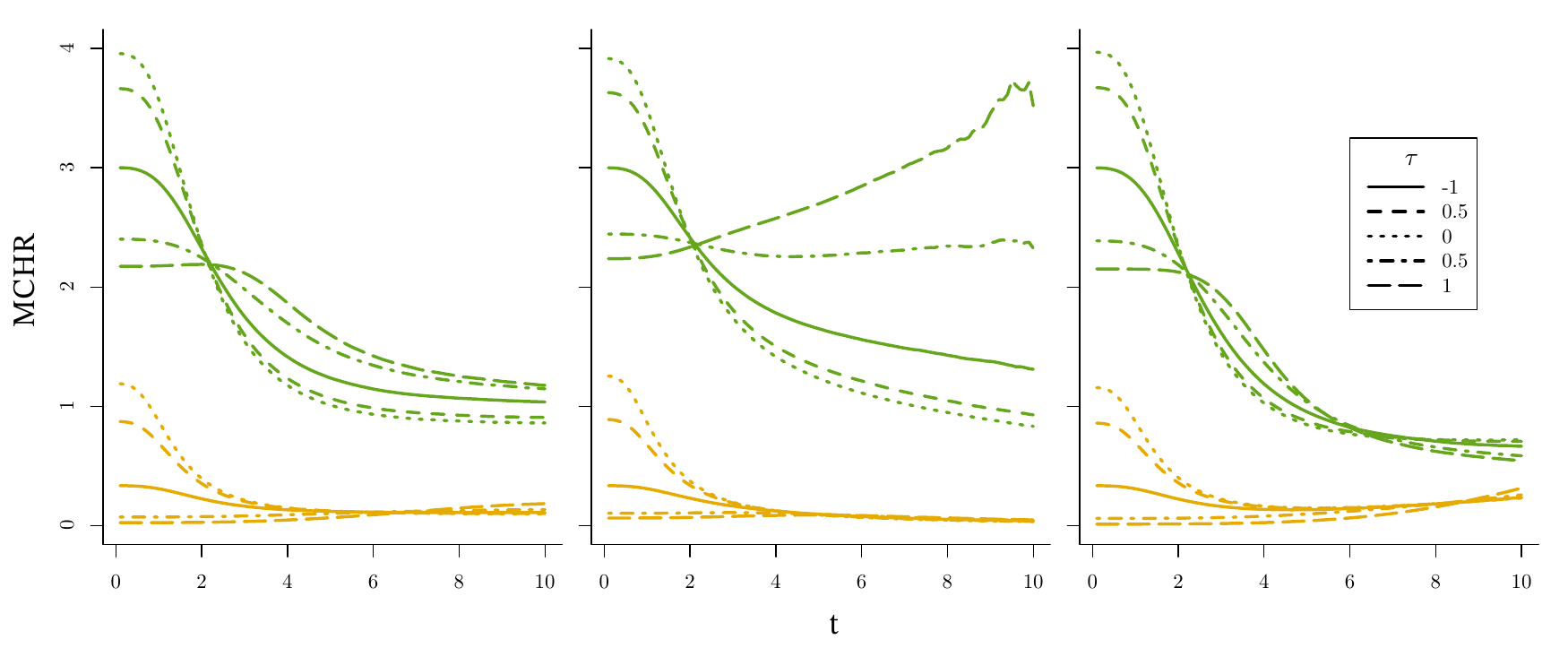}			
	\caption{\mchr  over time for $\lambda^{a}_{i}(t)=U_{0i}(U_{1i})^{a} \tfrac{t^{2}}{20}$, for a unit-variance gamma distributed $U_{1}$ with $\mathbb{E}[U_{1}]=3$ (green) or $\mathbb{E}[U_{1}]=\tfrac{1}{3}$ (orange), $U_{0}$ follows a gamma (left), inverse Gaussian (middle) or compound Poisson (right) distribution and where the joint distribution of $U_{0}$ and $U_{1}$ follows from a Gaussian copula with varying Kendall's $\tau$ correlation coefficient (see legend). }\label{CH5Fig5b}	
\end{figure}

 \noindent Note that for $\rho=0$, we recover the independent setting already shown in Figure \ref{CH5FrailDisc} that can be used for comparison. First of all, it is important to recall that when $U_{0}$ and $U_{1}$ are dependent the \chr equals $\mathbb{E}[U_{1}]+\text{cov}(U_{0},U_{1})$ as discussed in Section \ref{CH5sec:chr}. For \chrs greater than $1$, it becomes clear that the selection effect is more serious for cases with a high (positive) correlation between $U_{0}$ and $U_{1}$. The stronger selection effect is due to the higher variability of $U_{0}U_{1}$. For \chrs less than $1$, this trend is only true at short timescales, after which the frailty selection effect takes over since for this example a large fraction of the individuals, $p_{1}=0.9$, the effect is the same ($U_{1}=0.1$). When we use a continuous gamma distributed $U_{1}$ instead, the frailty selection effect is less apparent, as shown in Figure \ref{CH5Fig5b}.   So far, for \chrs larger than $1$, we have observed a monotonic \mchr. However, in the case of strong dependence between $U_{0}$ and $U_{1}$ ($\tau=-1$, and $\tau=1$ for inverse Gaussian frailty) $\mathbb{E}[U_{0}\mid T^{0}\geq t]$ decreases faster than $\mathbb{E}[U_{0}U_{1}\mid T^{1}\geq t]$ resulting in a non-monotonic trend for the \mchr. For a Gamma distributed $U_{1}$, in the case of inverse Gaussian distributed $U_{0}$ with $\tau=1$, the \mchr even equals a monotonic increasing function over time as shown in Figure \ref{CH5Fig5b}.

In this section, we have applied Theorem \ref{CH5th2} to several examples to illustrate the deviation of the \mchr from the \chr. In summary, even when the \chr is constant, an \ohr from an RCT equal to $x$ (at time $t$) can occur for different \chr values when the $(U_{0}, U_{1})$ distribution is unknown as summarized in Table \ref{CH5tab:overview}. 

\begin{table}[H]
		\centering
	\small
	\tabcolsep=0.10cm
	\captionsetup{width=\textwidth}
	\caption{Assuming no confounding, an \ohr (at time $t$) equal to $x$ can occur for all values of a constant \chr as a result of selection of the frailty $(U_{0})$ or modifier factor $(U_{1})$ which might be dependent. \label{CH5tab:overview}}
		\begin{tabular}{c||c|c|c} 
			&                  & Cause                    & Presented examples                                 \\ \hline   
			\multirow{3}{*}{$x>1$} & $\text{CHR}>x$   & Frailty or modifier selection  &                                                      \\
			& $1<\text{CHR}<x$ & Dependence $U_{0}$ and $U_{1}$ & inverse Gaussian frailty - $\tau=1$ (Figures \ref{CH5Fig5}, \ref{CH5Fig5b}) \\
			& $\text{CHR}<1$   & Frailty selection              & compound Poisson frailty                             \\ \hline   
			\multirow{3}{*}{$x<1$} & $\text{CHR}>1$   & Frailty or modifier selection  &                                                      \\
			& $x<\text{CHR}<1$ & Modifier selection             & gamma distributed modifier (Figure \ref{CH5Fig5b})                   \\
			& $\text{CHR}<x$   & Frailty selection              &         \\                                                		\end{tabular}%
\end{table}

\section{Implications for the Cox model}
\begin{figure}[h!]
	\centering
	\captionsetup{width=\textwidth}
	\includegraphics[width=0.6\textwidth]{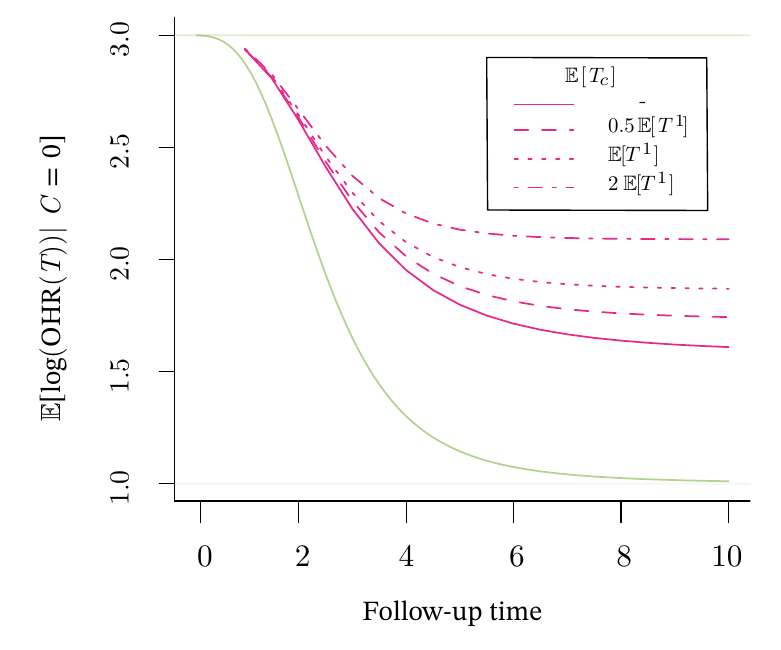} 
	\caption{Empirically obtained $\exp(\mathbb{E}[\log(\text{OHR}(T))\mid C=0])$ for increasing time-to follow up (solid pink), and for cases with an additional exponentially distributed censoring time ($T_{C}$, dashed pink), when $\lambda^{a}_{i}(t)=U_{0i}(U_{1i})^{a} \tfrac{t^{2}}{20}$,  $U_{1}\sim \text{BHN}(0.05,0.5,0.82,3.5)$ and $U_{0}\sim \Gamma(1,1)$. The \mchr  is also presented (green). 
		}
	\label{CH5figCox}	
\end{figure}

We have demonstrated that in the presence of frailty and effect heterogeneity, even when the \chr is constant, the \mchr varies over time. Then, the proportional hazards assumption will not hold for an \ohr from an RCT (that is in the absence of informative censoring equal to the \mchr as discussed at the start of Section \ref{CH5OHR}).   Despite the many options to deal with non-proportional hazards (see, e.g.~ \citet[Section~6.5]{Therneau2000} or \citet{Hess1994, Bennett1983, Wei2008}), in the majority of epidemiological time-to-event studies, the traditional Cox's proportional hazard model (that is thus misspecified) is fitted. The logarithm of the Cox estimate can be interpreted as the logarithm of the \ohr marginalized over the observed death times \citep{Schemper2009}, i.e.~$\mathbb{E}[\log(\text{OHR}(T))\mid C=0]$ for censoring indicator $C$.  The logarithm of the Cox estimate obtained from an RCT thus equals a time-weighted average of the logarithm of the \ohr. In the case of non-proportional hazards, even for uninformative censoring, the estimate is well-known to be affected by the censoring distribution. It differs from the average log hazard ratio $\mathbb{E}\left[\log(\text{OHR}(T))\right]$ 
\citep{Schemper2009, Xu2000, Boyd2012}. Therefore, the bias of the Cox estimate, when the estimand is the \chr, will depend on the joint distribution of $(U_0, U_1)$ as well as the censoring distribution. In most cases considered in Section \ref{CH5OHR}, the deviation of the \ohr from the \chr increased over time. For uninformative censoring, the probability of censoring increases over time, so the Cox estimate is closer to the \ohr at short times. In Figure \ref{CH5figCox}, this is demonstrated for the gamma-frailty case  (Figure \ref{CH5FrailDisc}, $\text{var}(U_{0}=1)$) by presenting empirically obtained  $\mathbb{E}[\log(\text{OHR}(T))\mid C=0]$ (with $1,000,000$ replications) 
based on a varying follow-up time  and loss to follow-up modelled with an exponential censoring-time distribution with different parameters.  The built-in selection bias of the \ohr (assuming no confounding) emphasizes the importance of testing the proportional hazard assumption. Serious bias will result in a time-varying \ohr, so that violation of the proportional hazard assumption can be verified when fitting a Cox model. However, when the actual \chr would be time-varying, the \ohr can be approximately constant when the selection effect and the change in \chr roughly cancels out (see, e.g.~\citet{Stensrud2019}). The data cannot be used to distinguish the latter case from the case where there is no heterogeneity and, thus, no selection effect. Similarly, as already mentioned at the end of Section \ref{CH5sec:chr}, when the \ohr would vary over time, we can never conclude whether this is the result of a time-varying causal effect or due to selection. However, the proportional hazard assumption would be violated in this case, and a standard Cox model is inappropriate.

\section{Discussion and concluding remarks}
In this paper, we have formalized how heterogeneity leads to deviation of the \mchr (see Equation \eqref{CH5nhr}) from the \chr of interest (see Definition \ref{CH5def:CHR}) due to the selection of both the individual frailty factor $(U_{0})$ and the individual effect modifier $(U_{1})$. This work generalizes frailty examples presented in the literature \citep{Hernan2010, Aalen2015, Stensrud2017}, by considering the possibility of multiplicative effect (on the hazard) heterogeneity that also results in non-exchangeability of exposed and unexposed individuals over time. As a result of the individual effect modifier ($U_{1}$), the individuals that survive in the exposed groups are expected to benefit more (or suffer less) from the exposure. At the same time $U_{0}\mid T^{1}\geq t$ will have a different distribution than $U_{0}\mid T^{0}\geq t$. When the \chr is larger than $1$ and $U_{0} \independent U_{1}$, the selection effects act in the same direction. On the other hand, when the \chr is smaller than $1$ and $U_{0} \independent U_{1}$, the selection effects can act in opposite directions so that the \mchr might be closer to the \chr than in the case without effect heterogeneity (see Figure \ref{CH5FrailDisc}). 

For data from an RCT, in the absence of informative censoring, the \ohr equals the studied \mchr so that all results directly relate to the \ohr. For observational data, the \ohr does not equal the \mchr due to confounding. However, when all confounders $\boldsymbol{L}$ are observed, i.e.~$T^{a} \independent A \mid \boldsymbol{L}$, one can study the conditional (on $\boldsymbol{L}$) \ohr that in turn is equal to the conditional \mchr. The presented theorems are valid while conditioning on $\boldsymbol{L}$.  

The intuition explained by \citet{Hernan2010} suggests that the \mchr underestimates the effect size, while the sign of the logarithms of the \mchr and the \chr are equal. However, we have shown that in the presence of effect heterogeneity, an \mchr equal to $x$ can occur both under $\text{CHR}>1$ as well as $\text{CHR}<1$ as summarized in Table \ref{CH5tab:overview}. Therefore, \ohrs from RCTs are not guaranteed to present a lower bound for the causal effect without making untestable assumptions on the $(U_{0}, U_{1})$ distribution. We have derived how the \mchr will evolve due to the selection of frailty and effect modifiers in Theorem \ref{CH5th2}. However, in practice, only the evolution of the \ohr can be found (assuming sufficient data is available). Even after assuming absence of confounding (e.g.~in case of a RCT) the \chr is non-identifiable without making (untestable) assumptions on the $(U_{0},U_{1})$ distribution as discussed at the end of Section \ref{CH5sec:chr}. 

We can thus not distinguish between a time-varying \chr without selection of $U_{0}$ and $U_{1}$ or a constant \chr with selection, see e.g \citet{Stensrud2020}. However, adjusting for other risk factors can lower the remaining variability of $U_{0}$ and $U_{1}$ so that the bias is reduced. Even for an RCT, it may thus help to focus on adjusted hazard ratios despite the absence of confounding. Nevertheless, adjusting for other risk factors will require more data and modelling decisions. We want to remark that although additive hazard models (when well-specified) do not suffer from the frailty selection as shown by \citet{Aalen2015}, these models will suffer from latent modifier selection in the presence of effect heterogeneity ($\mathbb{E}[U_{1} \mid T^{1}\geq t]>\mathbb{E}[U_{1}]$) as demonstrated in our companion paper \citet{Post2022c}. 

We hope that the discussed effect heterogeneity and formalization of the built-in selection bias of the \ohr (in case of an RCT) show the need to use more appropriate estimands. As suggested by others,  contrasts of the survival probabilities, the median, or the restricted mean survival time respectively of potential outcomes, are appropriate measures to quantify causal effects on time-to-event outcomes \citep{Hernan2010, Bartlett2020, Stensrud2019, Young2020}. 

		\bibliographystyle{spbasic}
		\bibliography{CHR}

\begin{thebibliography}{35}
\providecommand{\natexlab}[1]{#1}
\providecommand{\url}[1]{{#1}}
\providecommand{\urlprefix}{URL }
\expandafter\ifx\csname urlstyle\endcsname\relax
  \providecommand{\doi}[1]{DOI~\discretionary{}{}{}#1}\else
  \providecommand{\doi}{DOI~\discretionary{}{}{}\begingroup
  \urlstyle{rm}\Url}\fi
\providecommand{\eprint}[2][]{\url{#2}}

\bibitem[{Aalen et~al.(2015)Aalen, Cook, and R{\o}ysland}]{Aalen2015}
Aalen OO, Cook RJ, R{\o}ysland K (2015) {Does {C}ox analysis of a randomized
  survival study yield a causal treatment effect?} Lifetime Data Anal
  21(4):579--593

\bibitem[{Balan and Putter(2020)}]{Balan2020}
Balan TA, Putter H (2020) A tutorial on frailty models. Statistical Methods in
  Medical Research 29(11):3424--3454

\bibitem[{Bartlett et~al.(2020)Bartlett, Morris, Stensrud, Daniel,
  Vansteelandt, and Burman}]{Bartlett2020}
Bartlett JW, Morris TP, Stensrud MJ, Daniel RM, Vansteelandt SK, Burman CF
  (2020) The hazards of period specific and weighted hazard ratios. Statistics
  in Biopharmaceutical Research 12(4):518--519

\bibitem[{Bennett(1983)}]{Bennett1983}
Bennett S (1983) Analysis of survival data by the proportional odds model.
  Statistics in Medicine 2(2):273--277

\bibitem[{Bongers et~al.(2021)Bongers, Forr{\'e}, Peters, and
  Mooij}]{Bongers2021}
Bongers S, Forr{\'e} P, Peters J, Mooij JM (2021) {Foundations of structural
  causal models with cycles and latent variables}. The Annals of Statistics
  49(5):2885 -- 2915

\bibitem[{Boyd et~al.(2012)Boyd, Kittelson, and Gillen}]{Boyd2012}
Boyd AP, Kittelson JM, Gillen DL (2012) Estimation of treatment effect under
  non-proportional hazards and conditionally independent censoring. Statistics
  in Medicine 31(28):3504--3515

\bibitem[{Cox(1972)}]{Cox1972}
Cox DR (1972) Regression models and life-tables. Journal of the Royal
  Statistical Society Series B 34(2):187--220

\bibitem[{Daniel et~al.(2021)Daniel, Zhang, and Farewell}]{Daniel2020}
Daniel R, Zhang J, Farewell D (2021) Making apples from oranges: Comparing
  noncollapsible effect estimators and their standard errors after adjustment
  for different covariate sets. Biometrical Journal 63(3):528--557

\bibitem[{De~Neve and Gerds(2020)}]{DeNeve2020}
De~Neve J, Gerds TA (2020) On the interpretation of the hazard ratio in {C}ox
  regression. Biometrical Journal 62(3):742--750

\bibitem[{Didelez and Stensrud(2021)}]{Didelez2021}
Didelez V, Stensrud MJ (2021) On the logic of collapsibility for causal effect
  measures. Biometrical Journal

\bibitem[{Hern\'{a}n(2010)}]{Hernan2010}
Hern\'{a}n MA (2010) The hazards of hazard ratios. Epidemiology 21(1):13--15

\bibitem[{Hern{\'{a}}n and Robins(2020)}]{Hernan2019}
Hern{\'{a}}n MA, Robins JM (2020) {Causal Inference: What If.} Boca Raton:
  Chapman {\&} Hall/CRC

\bibitem[{Hern\'{a}n et~al.(2004)Hern\'{a}n, Hern\'{a}ndez-D\'{i}az, and
  Robins}]{Hernan2004a}
Hern\'{a}n MA, Hern\'{a}ndez-D\'{i}az S, Robins JM (2004) A structural approach
  to selection bias. Epidemiology 15(5):615--625

\bibitem[{Hern\'{a}n et~al.(2005)Hern\'{a}n, Cole, Margolick, Cohen, and
  Robins}]{Hernan2005b}
Hern\'{a}n MA, Cole SR, Margolick J, Cohen M, Robins JM (2005) Structural
  accelerated failure time models for survival analysis in studies with
  time-varying treatments. Pharmacoepidemiology and Drug Safety 14(7):477--491

\bibitem[{Hess(1994)}]{Hess1994}
Hess KR (1994) Assessing time-by-covariate interactions in proportional hazards
  regression models using cubic spline functions. Statistics in Medicine
  13(10):1045--1062

\bibitem[{Howe et~al.(2016)Howe, Cole, Lau, Napravnik, and Eron}]{Howe2016}
Howe CJ, Cole SR, Lau B, Napravnik S, Eron JJ (2016) Selection bias due to loss
  to follow up in cohort studies. Epidemiology 27(1):91--97

\bibitem[{Martinussen and Vansteelandt(2013)}]{Martinussen2013}
Martinussen T, Vansteelandt S (2013) On collapsibility and confounding bias in
  {C}ox and {A}alen regression models. Lifetime Data Analysis 19(3):279--296

\bibitem[{Martinussen et~al.(2020)Martinussen, Vansteelandt, and
  Andersen}]{Martinussen2020}
Martinussen T, Vansteelandt S, Andersen P (2020) Subtleties in the
  interpretation of hazard contrasts. Lifetime Data Analysis 26(4):833--855

\bibitem[{Nelsen(2006)}]{Nelsen2006}
Nelsen RB (2006) An Introduction to Copulas, 2nd edn. Springer

\bibitem[{Neyman(1990)}]{Neyman1990}
Neyman J (1990) {On the Application of Probability Theory to Agricultural
  Experiments. Essay on Principles. Section 9}. Statistical Science
  5(4):465--472

\bibitem[{Pearl(2009)}]{Pearl2009book}
Pearl J (2009) {Causality: Models, reasoning, and inference}, 2nd edn.
  Cambridge University Press, Cambridge

\bibitem[{Peters et~al.(2018)Peters, Janzing, and Sch{\"o}lkopf}]{Peters2017}
Peters J, Janzing D, Sch{\"o}lkopf B (2018) {Elements of causal inference:
  foundations and learning algorithms}. The MIT Press, Cambridge, Massachusetts

\bibitem[{Post et~al.(2022)Post, {van den Heuvel}, and Putter}]{Post2022c}
Post RAJ, {van den Heuvel} ER, Putter H (2022) Bias of the additive hazard
  model in the presence of causal effect heterogeneity

\bibitem[{Richardson and Robins(2013)}]{Richardson2013}
author (2013) {Single World Intervention Graphs (SWIGs): A Unification of the
  Counterfactual and Graphical Approaches to Causality}. Tech. Rep. 128,
  University of Washington

\bibitem[{Rubin(1974)}]{Rubin1974}
Rubin DB (1974) {Estimating causal effects of treatments in randomized and
  nonrandomized studies}. Journal of Educational Psychology 66(5):688--701

\bibitem[{Schemper et~al.(2009)Schemper, Wakounig, and Heinze}]{Schemper2009}
Schemper M, Wakounig S, Heinze G (2009) The estimation of average hazard ratios
  by weighted cox regression. Statistics in Medicine 28(19):2473--2489

\bibitem[{Sjölander et~al.(2016)Sjölander, Dahlqwist, and
  Zetterqvist}]{Sjolander2016}
Sjölander A, Dahlqwist E, Zetterqvist J (2016) A note on the noncollapsibility
  of rate differences and rate ratios. Epidemiology 27(3):356--359

\bibitem[{Sklar(1959)}]{Sklar1959}
Sklar A (1959) Fonctions de r\'{e}partition \`{a} n dimensions et leurs marges.
  Publications de l'Institut de statistique de l'Universit\'{e} de Paris
  8:229--231

\bibitem[{Stensrud and Hernán(2020)}]{Stensrud2020}
Stensrud MJ, Hernán MA (2020) Why test for proportional hazards? JAMA
  323(14):1401--1402

\bibitem[{Stensrud et~al.(2017)Stensrud, Valberg, R{\o}ysland, and
  Aalen}]{Stensrud2017}
Stensrud MJ, Valberg M, R{\o}ysland K, Aalen OO (2017) {Exploring selection
  bias by causal frailty models: The magnitude matters}. Epidemiology
  28(3):379--386

\bibitem[{Stensrud et~al.(2018)Stensrud, Aalen, Aalen, and
  Valberg}]{Stensrud2019}
Stensrud MJ, Aalen JM, Aalen OO, Valberg M (2018) {Limitations of hazard ratios
  in clinical trials}. European Heart Journal 40(17):1378--1383

\bibitem[{Thernau and Grambsch(2000)}]{Therneau2000}
Thernau TM, Grambsch PM (2000) Modeling Survival Data: Extending the Cox Model,
  1st edn. Springer

\bibitem[{Wei and Schaubel(2008)}]{Wei2008}
Wei G, Schaubel DE (2008) Estimating cumulative treatment effects in the
  presence of nonproportional hazards. Biometrics 64(3):724--732

\bibitem[{Xu and O'Quigley(2000)}]{Xu2000}
Xu R, O'Quigley J (2000) {Estimating average regression effect under
  non-proportional hazards }. Biostatistics 1(4):423--439

\bibitem[{Young et~al.(2020)Young, Stensrud, {Tchetgen Tchetgen}, and
  Hern{\'{a}}n}]{Young2020}
Young JG, Stensrud MJ, {Tchetgen Tchetgen} EJ, Hern{\'{a}}n MA (2020) {A causal
  framework for classical statistical estimands in failure-time settings with
  competing events}. Stat Med 39(8):1199--1236

\end{thebibliography}
		
\appendix

\section{Proofs}\label{CH5CH5:proofs}
\subsection{Proof of Theorem \ref{CH5th2}}
\begin{proof}
	\begin{align*}
	\frac{\mathbb{E}\left[\lambda_{i}^{a}(t)\right]}{\mathbb{E}\left[\lambda_{i}^{0}(t)\right]}&=\frac{\mathbb{E}\left[f_{0}(t,U_{0i})f_{1}(a,U_{1i},t)\right]}{\mathbb{E}\left[f_{0}(t,U_{0i})\right]}\\
	&= \frac{\int\lim_{h\rightarrow 0}h^{-1}\mathbb{P}\left(T^{a} \in [t,t+h) \mid T^{a} \geq t, U_{0}, U_{1} \right)dF_{U_{0}, U_{1i}}}{\int\lim_{h\rightarrow 0}h^{-1}\mathbb{P}\left(T^{0} \in [t,t+h) \mid  T^{0}\geq t, U_{0}  \right)dF_{U_{0}}}\\
	&= \frac{\int\lim_{h\rightarrow 0}h^{-1}\mathbb{P}\left(T^{a} \in [t,t+h) \mid T^{a} \geq t, U_{0}, U_{1},A=a \right)dF_{U_{0}, U_{1}}}{\int\lim_{h\rightarrow 0}h^{-1}\mathbb{P}\left(T^{0} \in [t,t+h) \mid  T^{0}\geq t, U_{0}, A=0  \right)dF_{U_{0}}}\\
	&= \frac{\int\lim_{h\rightarrow 0}h^{-1}\mathbb{P}\left(T \in [t,t+h) \mid T\geq t, U_{0}, U_{1},A=a \right)dF_{U_{0}, U_{1}}}{\int\lim_{h\rightarrow 0}h^{-1}\mathbb{P}\left(T \in [t,t+h) \mid  T\geq t, U_{0}, A=0  \right)dF_{U_{0}}}  
	\end{align*}
\end{proof}

\subsection{Proof of Theorem \ref{CH5th42}}
\begin{proof}
	By the law of total probability, 
	\begin{align*}
	&\frac{\lim_{h\rightarrow 0}h^{-1}\mathbb{P}\left(T^{a} \in [t,t+h) \mid T^{a}\geq t \right)}{\lim_{h\rightarrow 0}h^{-1}\mathbb{P}\left(T^{0} \in [t,t+h) \mid T^{0}\geq t \right)} =\\ 
	&\frac{\lim_{h\rightarrow 0}\int h^{-1}\mathbb{P}\left(T^{a} \in [t,t+h) \mid T^{a}\geq t, U_{0}, U_{1} \right)dF_{U_{0},U_{1}\mid T^{a}\geq t}}{\lim_{h\rightarrow 0}\int h^{-1}\mathbb{P}\left(T^{0} \in [t,t+h) \mid T^{0}\geq t, U_{0} \right) dF_{U_{0}\mid T^{0}\geq t}}
	\end{align*} First we focus on the integrand,  \begin{align*}
	&h^{-1}\mathbb{P}\left(T^{a} \in [t,t+h) \mid T^{a}\geq t, U_{0}, U_{1} \right)\\ &=  h^{-1}\frac{\mathbb{P}\left(T^{a} \geq  t \mid  U_{0}, U_{1} \right) - \mathbb{P}\left(T^{a} \geq  t + h \mid  U_{0}, U_{1} \right)}{\mathbb{P}\left(T^{a} \geq  t \mid  U_{0}, U_{1} \right)} \\ 
	&=  h^{-1}\left(1-\frac{\mathbb{P}\left(T^{a} \geq  t + h \mid  U_{0}, U_{1} \right)}{\mathbb{P}\left(T^{a} \geq  t \mid  U_{0}, U_{1} \right)}\right) \\ 
	&=  h^{-1}\left(1-\frac{\exp\left(- \int_{0}^{t+h} f_{0}(s,U_{0})f_{1}(s,U_{1},a) ds \right)}{\exp\left(- \int_{0}^{t} f_{0}(s,U_{0})f_{1}(s,U_{1},a) ds \right)}\right) \\
	&=  h^{-1}\left(1-\exp\left(- \int_{t}^{t+h} f_{0}(s,U_{0})f_{1}(s,U_{1},a) ds \right)\right) 
	\end{align*} For monotonic (increasing or decreasing) conditional hazard functions if $h_{2}<h_{1}$, then \begin{equation*}
	h_{1}^{-1}\left(1- \exp\left(- \int_{t}^{t+h_{1}} f_{0}(s,U_{0})f_{1}(s,U_{1},a) ds \right)\right) \leq  h_{2}^{-1}\left(1- \exp\left(- \int_{t}^{t+h_{2}} f_{0}(s,U_{0})f_{1}(s,U_{1},a) ds \right)\right)	\end{equation*} or 
	\begin{equation*}
	h_{1}^{-1}\left(1- \exp\left(- \int_{t}^{t+h_{1}} f_{0}(s,U_{0})f_{1}(s,U_{1},a) ds \right)\right) \geq  h_{2}^{-1}\left(1- \exp\left(- \int_{t}^{t+h_{2}} f_{0}(s,U_{0})f_{1}(s,U_{1},a) ds \right)\right)	\end{equation*} as the average integrated conditional hazard over the interval increases (or decreases). Then, the limit and integral can be interchanged by directly applying the monotone convergence theorem.
	
	For non-monotone conditional hazard functions, let $h \leq \tilde{{}h}$, and note 
	\begin{align*}
	h^{-1}\mathbb{P}\left(T^{a} \in [t,t+h) \mid T^{a}\geq t, U_{0}, U_{1} \right) &\leq h^{-1}\left(1- \exp\left(- h f_{0}(t^{*},U_{0})f_{1}(t^{*},U_{1},a)\right)\right)
	\end{align*} where
	\begin{equation*}
	\max_{s \in (t, t + \tilde{{}h})} f_{0}(s,U_{0})f_{1}(s,U_{1},a) = f_{0}(t^{*},U_{0})f_{1}(t^{*},U_{1},a).
	\end{equation*} Using the power series definition of the exponential function, \begin{align*}
	&h^{-1}\mathbb{P}\left(T^{a} \in [t,t+h) \mid T^{a}\geq t, U_{0}, U_{1} \right)\\
	&\leq h^{-1}\left(1- \frac{1}{\sum_{k=0}^{\infty}h^{k}  (f_{0}(t^{*},U_{0})f_{1}(t^{*},U_{1},a))^{k}\tfrac{1}{k!}}\right)\\
	&= h^{-1} \frac{\sum_{k=1}^{\infty}h^{k}  (f_{0}(t^{*},U_{0})f_{1}(t^{*},U_{1},a))^{k}\tfrac{1}{k!}}{\sum_{k=0}^{\infty}h^{k}  (f_{0}(t^{*},U_{0})f_{1}(t^{*},U_{1},a))^{k}\tfrac{1}{k!}}\\
	&= f_{0}(t^{*},U_{0})f_{1}(t^{*},U_{1},a) \frac{\sum_{k=1}^{\infty}h^{k-1}  (f_{0}(t^{*},U_{0})f_{1}(t^{*},U_{1},a))^{k-1}\tfrac{1}{k!}}{\sum_{k=0}^{\infty}h^{k}  (f_{0}(t^{*},U_{0})f_{1}(t^{*},U_{1},a))^{k}\tfrac{1}{k!}}\\
	&= f_{0}(t^{*},U_{0})f_{1}(t^{*},U_{1},a) \frac{\sum_{k=0}^{\infty}h^{k}  (f_{0}(t^{*},U_{0})f_{1}(t^{*},U_{1},a))^{k}\tfrac{1}{(k+1)!}}{\sum_{k=0}^{\infty}h^{k}  (f_{0}(t^{*},U_{0})f_{1}(t^{*},U_{1},a))^{k}\tfrac{1}{k!}}\\
	&< f_{0}(t^{*},U_{0})f_{1}(t^{*},U_{1},a).
	\end{align*} Furthermore, 
	$\mathbb{E}\left[ f_{0}(t^{*},U_{0})f_{1}(t^{*},U_{1},a) \mid T^{a}\geq t\right]<\infty$ when $\forall h \in (0,\tilde{{}h}):  \mathbb{E}[ f_{0}(t+h,U_{0}) f_{1}(t+h,U_{1},a) \mid T^{a}\geq t]<\infty$. Then we can change the order of the limit and integral by application of the dominated convergence theorem and conclude, 
	\begin{equation*}
	\frac{\lim_{h\rightarrow 0}h^{-1}\mathbb{P}\left(T^{a} \in [t,t+h) \mid T^{a}\geq t \right)}{\lim_{h\rightarrow 0}h^{-1}\mathbb{P}\left(T^{0} \in [t,t+h) \mid T^{0}\geq t \right)}
	= \frac{\mathbb{E}\left[f_{0}(t,U_{0})f_{1}(t,U_{1},a)\mid T^{a}\geq t\right]}{\mathbb{E}\left[f_{0}(t,U_{0})\mid T^{0}\geq t\right]}.
	\end{equation*}  
	
	Applying Bayes rule, we obtain \begin{align*}
	&\mathbb{E}\left[f_{0}(t,U_{0})f_{1}(t,U_{1},a)\mid T_{i}^{a}\geq t\right]\\
	&= \int f_{0}(t,U_{0})f_{1}(t,U_{1},a) dF_{(U_{0},U_{1})\mid T^{a}\geq t}\\
	&= \int f_{0}(t,U_{0})f_{1}(t,U_{1},a)\frac{ \mathbb{P}(T^a\geq t\mid U_{0},U_{1})}{\mathbb{P}(T^{a}\geq t)}dF_{(U_{0},U_{1})}\\
	&= \int f_{0}(t,U_{0})f_{1}(t,U_{1},a)\frac{\exp(-\int_{0}^{t}f_{0}(s,U_{0})f_{1}(s,U_{1},a)ds)}{\int \int\exp(-\int_{0}^{t}f_{0}(s,U_{0})f_{1}(s,U_{1},a)ds)dF_{(U_{0},U_{1})}}dF_{(U_{0},U_{1})}.
	\end{align*} Furthermore, \begin{align*}
	& \mathbb{E}\left[f_{0}(t,U_{0})\mid T^{a}\geq t\right]\\
        &=\int f_{0}(t,U_{0}) dF_{U_{0}\mid T_{i}^{a}\geq t}\\
	&=\int f_{0}(t,U_{0})\frac{ \mathbb{P}(T^a\geq t\mid U_{0})}{\mathbb{P}(T^{a}\geq t)}dF_{U_{0}}\\
	&= \int f_{0}(t,U_{0})\frac{\int\exp(-\int_{0}^{t}f_{0}(s,U_{0})f_{1}(s,U_{1},a)ds)dF_{U_{1}\mid U_{0}}}{ \int  \exp(-\int_{0}^{t}f_{0}(s,U_{0})f_{1}(s,U_{1},a)ds) dF_{(U_{0},U_{1})}}dF_{U_{0}},
	\end{align*} such that for $a=0$,
	\begin{equation*}
	\mathbb{E}\left[f_{0}(t,U_{0})\mid T^{0}\geq t\right]=\int f_{0}(t,U_{0}) \frac{\exp(-\int_{0}^{t}f_{0}(s,U_{0})ds)}{\int \exp(-\int_{0}^{t}f_{0}(s,U_{0})ds)dF_{U_{0}}}dF_{U_{0}}.
	\end{equation*} The ratio of both gives us the result. 
\end{proof}

\subsection{Proof of Corollary \ref{CH5cor411}}

\begin{proof}
	By Theorem \ref{CH5th42},
	\begin{equation*}
	\frac{\lim_{h\rightarrow 0}h^{-1}\mathbb{P}\left(T^{a} \in [t,t+h) \mid T^{a}\geq t \right)}{\lim_{h\rightarrow 0}h^{-1}\mathbb{P}\left(T^{0} \in [t,t+h) \mid T^{0}\geq t \right)}
	= \frac{\mathbb{E}\left[f_{0}(t,U_{0})f_{1}(t,U_{1},a)\mid T^{a}\geq t\right]}{\mathbb{E}\left[f_{0}(t,U_{0})\mid T^{0}\geq t\right]}.
	\end{equation*}  \noindent and equals
	\begin{equation*}
	\int f_{0}(t,U_{0})f_{1}(t,U_{1},a)\tfrac{\exp(-\Lambda^{a}(t,U_{0},U_{1})}{\int\int \exp(-\Lambda^{a}(t,U_{0},U_{1})dF_{(U_{0},U_{1})}}dF_{(U_{0},U_{1})}
	\left(\int f_{0}(t,U_{0})\tfrac{\exp(-\Lambda^{0}(t,U_{0}))}{\int \exp(-\Lambda^{0}(t,U_{0}))dF_{U_{0}}}dF_{U_{0}}\right)^{-1},
	\end{equation*} where $\Lambda^{a}(t,u_{0},u_{1})=\int_{0}^{t}f_{0}(s,u_{0})f_{1}(s,u_{1},a)ds$ and thus $\Lambda^{0}(t,u_{0})=\int_{0}^{t}f_{0}(s,u_{0})ds$. 
	\noindent As $U_{1}$ is degenerate, 
	\begin{align*}
	&\int f_{0}(t,U_{0})\tfrac{\exp(-\int_{0}^{t}f_{0}(s,U_{0})f_{1}(a,s)ds)} {\int \exp(-\int_{0}^{t}f_{0}(s,U_{0})f_{1}(a,s)ds)dF_{U_{0}}}dF_{U_{0}}
	\left(\int f_{0}(t,U_{0})\tfrac{\exp(-\int_{0}^{t}f_{0}(s,U_{0})ds)}{\int \exp(-\int_{0}^{t}f_{0}(s,U_{0}))dF_{U_{0}}}dF_{U_{0}}\right)^{-1}f_{1}(t,a)\\
	& = \frac{\mathbb{E}\left[f_{0}(t,U_{0})\mid T_{i}^{a}\geq t\right]}{\mathbb{E}\left[f_{0}(t,U_{0})\mid T_{i}^{0}\geq t\right]}f_{1}(t,a). 
	\end{align*}
\end{proof}

\subsection{Proof of Lemma \ref{CH5lemma1}}
\begin{proof}
	By Bayes rule, the probability density $f_{U_{0}\mid T^{a}\geq t}$ equals
	\begin{align*}
	f_{U_{0}\mid T^{a}\geq t}(u_{0}) &= \frac{\mathbb{P}(T^{a}\geq t\mid U_{0}=u_{0})f(u_{0})}{\int \mathbb{P}(T^{a}\geq t\mid U_{0})dF_{U_{0}}}\\
	&= \frac{\exp(-u_{0}\int_{0}^{t}\lambda_{0}(s)f_{1}(a,s) ds)f_{U_{0}}(u_{0})}{\int \exp(-U_{0}\int_{0}^{t}\lambda_{0}(s)f_{1}(a,s))dF_{U_{0}}}.
	\end{align*}
		\noindent Such that the Laplace transform of $U_{0}\mid T^{a}\geq t$ can be written as
		\begin{align*}
	\mathcal{L}_{U_{0}\mid T^{a}\geq t}(c) &= \mathbb{E}[\exp(-U_{0}c)\mid T^{a}\geq t]\\
	&= \int \exp(-u_{0}c)  f_{U_{0}\mid T^{a}\geq t}(u_{0}) du_{0}\\
	&= \int \exp(-u_{0}c) \frac{\exp(-u_{0}\int_{0}^{t}\lambda_{0}(s)f_{1}(a,s) ds)f_{U_{0}}(u_{0})}{\int \exp(-U_{0}\int_{0}^{t}\lambda_{0}(s)f_{1}(a,s))dF_{U_{0}}}\\
	&= \int  \frac{\exp(-u_{0}(c+\int_{0}^{t}\lambda_{0}(s)f_{1}(a,s) ds))f_{U_{0}}(u_{0})}{\int \exp(-U_{0}\int_{0}^{t}\lambda_{0}(s)f_{1}(a,s))dF_{U_{0}}}\\
	&=  \frac{\mathbb{E}\left[\exp(-U_{0}(c+\int_{0}^{t}\lambda_{0}(s)f_{1}(a,s) ds))\right]}{\mathbb{E}\left[ \exp(-U_{0}\int_{0}^{t}\lambda_{0}(s)f_{1}(a,s))\right]}\\
	&= \frac{\mathcal{L}_{U_{0}}(c+\int_{0}^{t}\lambda_{0}(s)f_{1}(a,s) ds)}{\mathcal{L}_{U_{0}}(\int_{0}^{t}\lambda_{0}(s)f_{1}(a,s) ds)}. 
	\end{align*}
	\noindent Since for a random variable $X$, $\mathbb{E}[X]= -\mathcal{L}_{X}^{'}(0)$, 
	\begin{equation*}
	\mathbb{E}[U_{0}\mid T^{a}\geq t] = -\mathcal{L}_{U_{0}\mid T^{a}\geq t}^{'}(0) = -\frac{\mathcal{L}_{U_{0}}^{'}(\int_{0}^{t}\lambda_{0}(s)f_{1}(a,s)ds)}{\mathcal{L}_{U_{0}}(\int_{0}^{t}\lambda_{0}(s)f_{1}(a,s)ds)}. 
	\end{equation*}
\end{proof}

\subsection{Proof of Corollary \ref{CH5cor421}}
\begin{proof}
	By Theorem \ref{CH5th42},
	\begin{equation*}
	\frac{\lim_{h\rightarrow 0}h^{-1}\mathbb{P}\left(T^{a} \in [t,t+h) \mid T^{a}\geq t \right)}{\lim_{h\rightarrow 0}h^{-1}\mathbb{P}\left(T^{0} \in [t,t+h) \mid T^{0}\geq t \right)}
	= \frac{\mathbb{E}\left[f_{0}(t,U_{0})f_{1}(t,U_{1})\mid T^{a}\geq t\right]}{\mathbb{E}\left[f_{0}(t,U_{0})\mid T^{0}\geq t\right]}.
	\end{equation*}  \noindent and equals
	\begin{equation*}
	\int f_{0}(t,U_{0})f_{1}(t,U_{1},a)\tfrac{\exp(-\Lambda^{a}(t,U_{0},U_{1})}{\int\int \exp(-\Lambda^{a}(t,U_{0},U_{1})dF_{(U_{0},U_{1})}}dF_{(U_{0},U_{1})}
	\left(\int f_{0}(t,U_{0})\tfrac{\exp(-\Lambda^{0}(t,U_{0}))}{\int \exp(-\Lambda^{0}(t,U_{0}))dF_{U_{0}}}dF_{U_{0}}\right)^{-1},
	\end{equation*} where $\Lambda^{a}(t,u_{0},u_{1})=\int_{0}^{t}f_{0}(s,u_{0})f_{1}(s,u_{1},a)ds$ and thus $\Lambda^{0}(t,u_{0})=\int_{0}^{t}f_{0}(s,u_{0})ds$. 
	\noindent As now $U_{0}$ is degenerate, 
	\begin{equation*}
	\int f_{1}(t,U_{1},a)\tfrac{\exp(-\Lambda^{a}(t,U_{1}))}{\int \exp(-\Lambda^{a}(t,U_{1}))dF_{U_{1}}} dF_{U_{1}} = \mathbb{E}\left[f_{1}(t,U_{1},a)\mid T_{i}^{a}\geq t\right],
	\end{equation*}
	where $\Lambda^{a}(t,u_{1})=\int_{0}^{t}\lambda_{0}(s)f_{1}(s,u_{1},a)ds$.
\end{proof}

\section{Laplace transforms}\label{CH5app:LT}
\subsection{Gamma}
If $X\sim \Gamma(k,\theta)$, then $\mathbb{E}[X]=k\theta$, $\text{var}[X]=k\theta^{2}$,
\begin{equation*}
\mathcal{L}_{X}(c)=(1+\theta c)^{-k},
\end{equation*}

\begin{equation*}
\mathcal{L}^{'}_{X}(c)= -\frac{\theta k}{\theta c +1}\mathcal{L}_{X}(c),
\end{equation*}
\noindent and
\begin{equation*}
\mathcal{L}^{''}_{X}(c)= \frac{\theta^{2}k(k+1)}{(\theta c +1)^{2}}\mathcal{L}_{X}(c).
\end{equation*} \noindent When $f_{\lambda}(t,U_{0i},U_{1i},a) = U_{0i}\lambda_{0}(t)f_{1}(t,a)$, by Lemma \ref{CH5lemma1}, 
\begin{equation}
\mathbb{E}\left[U_{0} \mid T^{a}\geq t\right]=
\frac{\theta  k}{\theta\Lambda^{a}(t) +1}, 
\end{equation} \noindent where $\Lambda^{a}(t)=\int_{0}^{t}\lambda_{0}(s)f_{1}(a,s)ds$. If $\mathbb{E}[U_{0}]=c$ and $\text{var}(U_{0})=\theta$, then $k=\frac{c^2}{\theta_{0}}$, $\theta=\frac{\theta_{0}}{c}$ and as such
\begin{equation}
\mathbb{E}\left[U_{0} \mid T^{a}\geq t\right]=
\frac{c}{\tfrac{\theta_{0}}{c}\Lambda^{a}(t) +1}. 
\end{equation} \noindent In particular, when $k=\theta^{-1}$, then $\mathbb{E}[U_{0}]=1$,  $\text{var}(U_{0})=\theta$ and 
\begin{equation}
\mathbb{E}\left[U_{0} \mid T^{a}\geq t\right]=
\left(\theta\Lambda^{a}(t) +1\right)^{-1}. 
\end{equation}

\subsection{Inverse Gaussian}
If $X\sim \text{IG}(\mu,\lambda)$, then $\mathbb{E}[X]=\mu$, $\text{var}[X]=\tfrac{\mu^3}{\lambda}$,
\begin{equation*}
\mathcal{L}_{X}(c)=\exp \left(\frac{\lambda}{\mu}\left(1-\sqrt{1+\frac{2\mu^{2}c}{\lambda}}\right)\right),
\end{equation*}
\begin{equation*}
\mathcal{L}^{'}_{X}(c)=-\frac{\mu}{\sqrt{\tfrac{2\mu^2c}{\lambda}+1}}\mathcal{L}_{X}(c), 
\end{equation*}
\noindent and 
\begin{equation*}
\mathcal{L}^{''}_{X}(c)=\left(\mu^{2} \frac{\lambda}{\lambda + 2 \mu^{2}c}+\frac{\mu}{\lambda(\tfrac{2\mu^{2}c}{\lambda}+1)^{\tfrac{3}{2}}}\right)\mathcal{L}_{X}(c).
\end{equation*} \noindent When $f_{\lambda}(t,U_{0i},U_{1i},a) = U_{0i}\lambda_{0}(t)f_{1}(t,a)$, by Lemma \ref{CH5lemma1}, \begin{equation}
\mathbb{E}\left[U_{0} \mid T^{a}\geq t\right]= \frac{\mu }{\sqrt{\frac{2 \Lambda^{a}(t) \mu ^2}{\lambda }+1}}, 
\end{equation} \noindent where $\Lambda^{a}(t)=\int_{0}^{t}\lambda_{0}(s)f_{1}(a,s)ds$. If $\text{var}(U_{0})=\theta_{0}$, then $\lambda=\frac{\mu^3}{\theta_{0}}$  and as such
\begin{equation}
\mathbb{E}\left[U_{0} \mid T^{a}\geq t\right]= \frac{\mu }{\sqrt{2 \Lambda^{a}(t) \tfrac{\theta_{0}}{\mu} +1}}.
\end{equation} \noindent In particular when $\mu=1$ and $\lambda=\theta_{0}^{-1}$, then $\mathbb{E}[U_{0}]=1$, $\text{var}(U_{0})=\theta_{0}$ and 
\begin{equation}
\mathbb{E}\left[U_{0} \mid T^{a}\geq t\right]= \left(2 \theta_{0} \Lambda^{a}(t)+1\right)^{-\tfrac{1}{2}}. 
\end{equation}

\subsection{Compound Poisson}
If $X\sim \text{CPois}(\rho, \eta, \nu)$, then $X= \sum_{i=1}^{N}Y_{i}$, where $N\sim Poi(\rho)$, and $Y \sim \Gamma(\eta,\nu)$,  $\mathbb{E}[X]=\rho\eta\nu$, $\text{var}[X]=\rho\eta\nu^{2}+\eta^2\nu^2\rho$,
\begin{equation*}
\mathcal{L}_{X}(c)=\exp \left(\rho\left(\left(\frac{\nu^{-1}}{\nu^{-1}+c}\right)^{\eta}-1\right)\right),
\end{equation*}
\begin{equation*}
\mathcal{L}^{'}_{X}(c)=-\rho\eta\nu  \left(\frac{1}{c \nu +1}\right)^{\eta +1}\mathcal{L}_{X}(c),
\end{equation*}
\noindent and
\begin{equation*}
\mathcal{L}^{''}_{X}(c)=\eta  \nu ^2 \rho  \left(\frac{1}{c \nu +1}\right)^{\eta +2} \left(\eta  \rho  \left(\frac{1}{c \nu +1}\right)^{\eta }+\eta +1\right)\mathcal{L}_{X}(c).
\end{equation*}
\noindent When $f_{\lambda}(t,U_{0i},U_{1i},a) = U_{0i}\lambda_{0}(t)f_{1}(t,a)$, by Lemma \ref{CH5lemma1}, 
\begin{equation}
\mathbb{E}\left[U_{0} \mid T^{a}\geq t\right]=
\rho\eta\nu  \left(\frac{1}{\Lambda^{a}(t) \nu +1}\right)^{\eta +1}, 
\end{equation} \noindent where $\Lambda^{a}(t)=\int_{0}^{t}\lambda_{0}(s)f_{1}(a,s)ds$.  Let  $\mathbb{E}[U_{0}]=c$ and  $\text{var}(U_{0})=\theta_{0}$, then $\rho=\frac{c^{2}(1+\eta)}{\eta \theta_{0}}$, $\nu=\frac{\theta_{0}}{c(1+\eta)}$  and as such
\begin{equation}
\mathbb{E}\left[U_{0} \mid T^{a}\geq t\right]=
c  \left(\frac{1}{\Lambda^{a}(t) \tfrac{\theta_{0}}{c(1+\eta)} +1}\right)^{\eta +1}. 
\end{equation} \noindent In particular, when $\eta=\tfrac{1}{2}$, $\rho=3\theta^{-1}$  and $\nu=\tfrac{2}{3} \theta,$ then $\mathbb{E}[U_{0}]=1$,  $\text{var}(U_{0})=\theta$ and 
\begin{equation}
\mathbb{E}\left[U_{0} \mid T^{a}\geq t\right]=
(1+\Lambda^{a}(t)\tfrac{2}{3}\theta)^{-\tfrac{3}{2}}. 
\end{equation} 

\subsection{Discrete}\label{CH5discrete}
Let $\mathbb{P}(X=\mu_{i})=p_{i}$ for $i>0, i\leq k$, 
\begin{equation*}
\mathcal{L}_{X}(c)=\sum_{i=1}^{n} p_{i} \exp \left(-c \mu_{i} \right),
\end{equation*}
\begin{equation*}
\mathcal{L}^{'}_{X}(c)=\sum_{i=1}^{n} -\mu_{i}p_{i} \exp \left(-c \mu_{i} \right),
\end{equation*} \noindent and
\begin{equation*}
\mathcal{L}^{''}_{X}(c)=\sum_{i=1}^{n} \mu_{i}^{2} p_{i} \exp \left(-c \mu_{i} \right).
\end{equation*}
\noindent Furthermore, if $Y\sim F$, then
\begin{equation*}
\mathcal{L}_{XY}(c)= \sum_{i=1}^{n} p_{i} \mathcal{L}_{Y}(c\mu_{i}),
\end{equation*}
\begin{equation*}
\mathcal{L}_{XY}^{'}(c)= \sum_{i=1}^{n} p_{i}\mu_{i} \mathcal{L}_{Y}^{'}(c\mu_{i}),
\end{equation*} \noindent and
\begin{equation*}
\mathcal{L}_{XY}^{''}(c)= \sum_{i=1}^{n} p_{i}\mu_{i}^{2} \mathcal{L}_{Y}^{''}(c\mu_{i}).
\end{equation*} 
\noindent Let $U_{1}$ equal $\mu_{1}$, $\mu_{2}$ or $1$ with probability $p_{1}$, $p_{2}$ and $1-p_{1}-p_{2}$ respectively. If 
\begin{equation*}
p_{2}=\frac{\left(\mu-\mu _1 p_1+p_1-1\right){}^2}{\mu^2-2 \mu-\mu _1^2 p_1+2 \mu _1 p_1-p_1+\theta_{1}+1},
\end{equation*} \noindent and
\begin{equation*}
\mu_{2}= \frac{\mu^2-\mu-\mu _1^2 p_1+\mu _1 p_1+\theta_{1}}{\mu-\mu _1 p_1+p_1-1},
\end{equation*} \noindent such that $p_{2} \in [0,1]$ and $\mu_{2}\geq 1$, then $\mathbb{E}[U_{1}]=\mu$ and $\text{var}(U_{1})=\theta_{1}$. 
When, $f_{\lambda}(t,U_{0i},U_{1i},a) = \lambda_{0}(t)(U_{1i})^{a}f_{1}(t,a)$, by Lemma \ref{CH5lemma1},  
\begin{equation}
\mathbb{E}\left[U_{1} \mid T^{1}\geq t\right]=\frac{\mu _1 p_1 e^{-\Lambda^{1}(t) \mu _1}+\mu _2 p_2 e^{-\Lambda^{1}(t) \mu _2}+\left(1-p_1-p_2\right)e^{-\Lambda^{1}(t)} }{p_1 e^{-\Lambda^{1}(t) \mu _1}+p_2 e^{-\Lambda^{1}(t) \mu _2}+\left(1-p_1-p_2\right)e^{-\Lambda^{1}(t)} }, 
\end{equation} \noindent where $\Lambda^{1} = \int_{0}^{t} \lambda_{0}(s) f_{1}(1,s) ds$. 

Furthermore, when $f_{\lambda}(t,U_{0i},U_{1i},a) = U_{0i} \lambda_{0}(t) (U_{1i})^{a} f_{1}(t,a)$, and $U_{0}$ and $U_{1}$ are independent, by Lemma \ref{CH5lemma1}, 
\begin{equation}
\mathbb{E}\left[U_{0}U_{1} \mid T^{1}\geq t\right] = \frac{p_{1}\mu_{1}\mathcal{L}_{U_{0}}^{'}(\Lambda^{1}\mu_{1})+
	p_{2}\mu_{2}\mathcal{L}_{U_{0}}^{'}(\Lambda^{1}\mu_{2})+
	(1-p_{1}-p_{2})\mathcal{L}_{U_{0}}^{'}(\Lambda^{1})}{p_{1}\mathcal{L}_{U_{0}}(\Lambda^{1}\mu_{1})+
	p_{2}\mathcal{L}_{U_{0}}(\Lambda^{1}\mu_{2})+
	(1-p_{1}-p_{2})\mathcal{L}_{U_{0}}(\Lambda^{1})}.
\end{equation} 

\section{Survival curves examples}\label{CH5app:surv}
The survival curves of $T^{0}$ and $T^{1}$ for the examples presented in this paper can be expressed in terms of the Laplace transforms presented in the previous section as shown in Table \ref{CH5tabsurv}, where $\Lambda_{0}(t)=\tfrac{t^{3}}{60}$. For the example where $U_{0} \not\independent U_{1}$, the survival curve for $T^{1}$ is obtained empirically from simulation. For data from a RCT $T^{a} \overset{d}{=} T\mid A=a$. 

\begin{table}[h!]
		\centering
	\small
	\tabcolsep=0.10cm
	\caption{Survival curves for $T^{0}$ and $T^{1}$ for different $\lambda_{i}^{a}(t)$ where $\Lambda_{0}(t)=\tfrac{t^{3}}{60}$.}\label{CH5tabsurv}
	\begin{tabular}{c||c|c} 
		$\lambda_{i}^{a}(t)$ & $S_{T^{0}}(t)$                        & $S_{T^{1}}(t)$                         \\ \hline
		$U_{0i}\tfrac{t^2}{20}c^{a}$               & $\mathcal{L}_{U_{0}}(\Lambda_{0}(t))$ & $\mathcal{L}_{U_{0}}(c\Lambda_{0}(t))$ \\ \hline
		$\tfrac{t^2}{20}(U_{1i})^{a}$              & $\exp(-\Lambda_{0}(t))$  &  $\mathcal{L}_{U_{1}}(\Lambda_{0}(t))$              \\ \hline
		$U_{0i}\tfrac{t^2}{20}(U_{1i})^{a}$                 & $\mathcal{L}_{U_{0}}(\Lambda_{0}(t))$ & $\mathcal{L}_{U_{0}U_{1}}(\Lambda_{0}(t))$                                 
	\end{tabular}%
	
\end{table}

\begin{figure}[h!]
	\centering
	\captionsetup{width=0\textwidth}
	\includegraphics[width=0.75\textwidth]{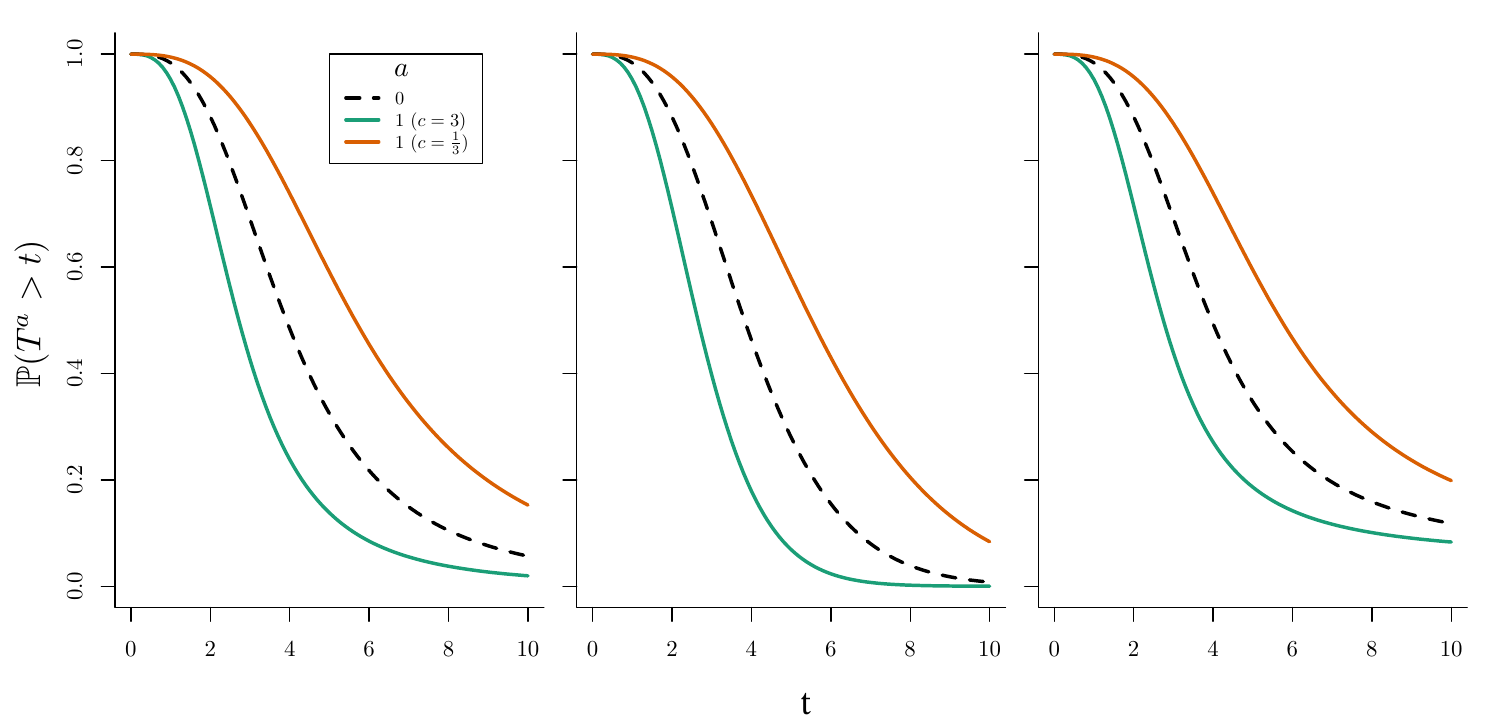} 
	\caption{Survival curves for $T^{0}$ (black), and $T^{1}$ (colored) when $\lambda^{a}_{i}(t)=U_{0i} \tfrac{t^{2}}{20}c^{a}$ for $c=3$ (green) and $c=\tfrac{1}{3}$ (orange), when $U_{0}$ follows a gamma, inverse Gaussian or compound Poisson distribution (from left to right) with variance $1$. The corresponding \mchr curves were presented in Figure \ref{CH5Frail}.  }\label{CH5Surv1}
\end{figure}
\begin{figure}[h!]
	\centering
	\captionsetup{width=\textwidth}
	\includegraphics[width=\textwidth]{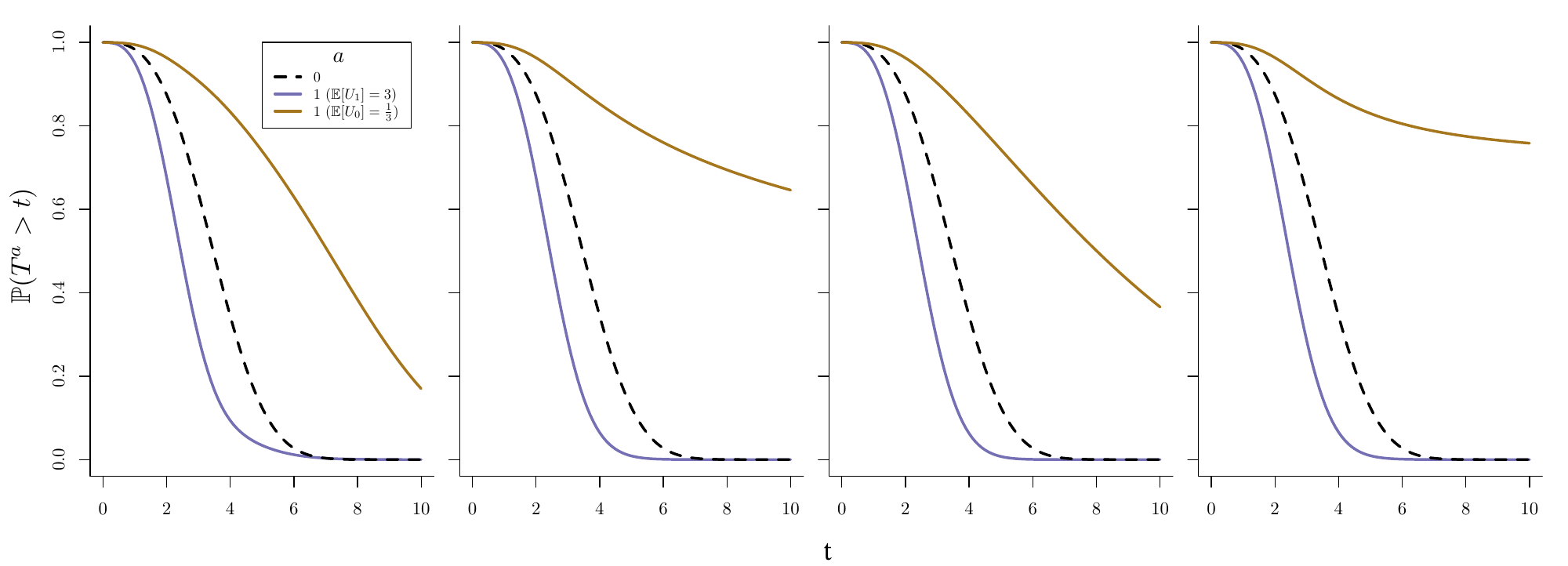} 
	\caption{
		Survival curves for $T^{0}$ (black), and $T^{1}$ (colored) when $\lambda^{a}_{i}(t)=U_{0i} \tfrac{t^{2}}{20}(U_{1i})^{a}$, when $U_{1}$ follows a BHN, gamma, inverse Gaussian or compound Poisson distribution (from left to right) with expectation $3$ (blue) or $\tfrac{1}{3}$ (brown) and variance  $1$. The corresponding \mchr curves were presented in Figure \ref{CH5Het}. }\label{CH5Surv2}
\end{figure}
\begin{figure}
	\centering
	\captionsetup{width=\textwidth}
	\includegraphics[width=0.75\textwidth]{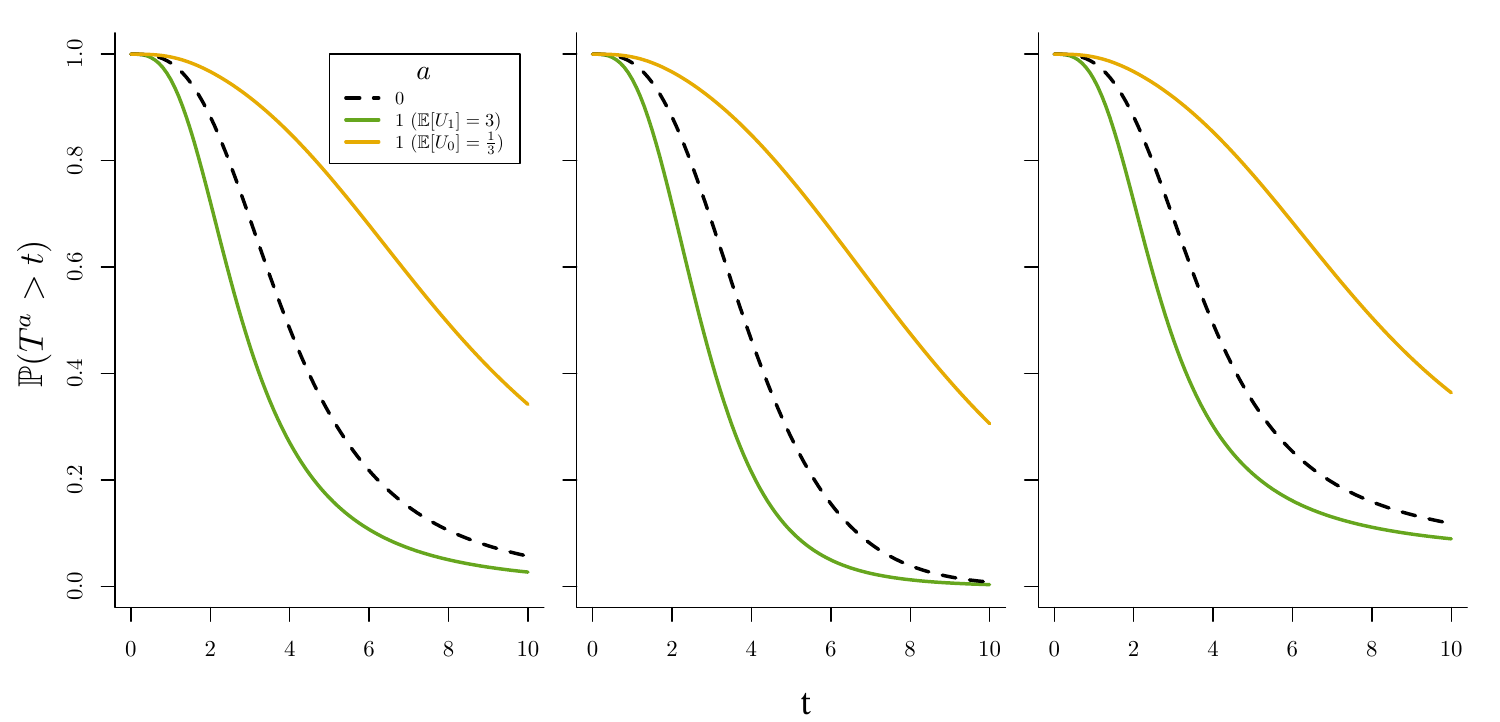} 
	\caption{Survival curves for $T^{0}$ (black), and $T^{1}$ (colored) when $\lambda^{a}(t)_{i}=U_{0i}(U_{1i})^{a} \tfrac{t^{2}}{20}$ for a unit-variance BHN distributed $U_{1}$ with $\mathbb{E}[U_{1}]=3$ (green) and $\mathbb{E}[U_{1}]=\tfrac{1}{3}$ (orange), when $U_{0}$ follows a gamma, inverse Gaussian or compound Poisson distribution (from left to right) with variance $1$. The corresponding \mchr curves were presented in Figure \ref{CH5FrailDisc}.}\label{CH5Surv3}
\end{figure} 
\begin{figure}
	\centering
	\captionsetup{width=\textwidth}
	\includegraphics[width=0.75\textwidth]{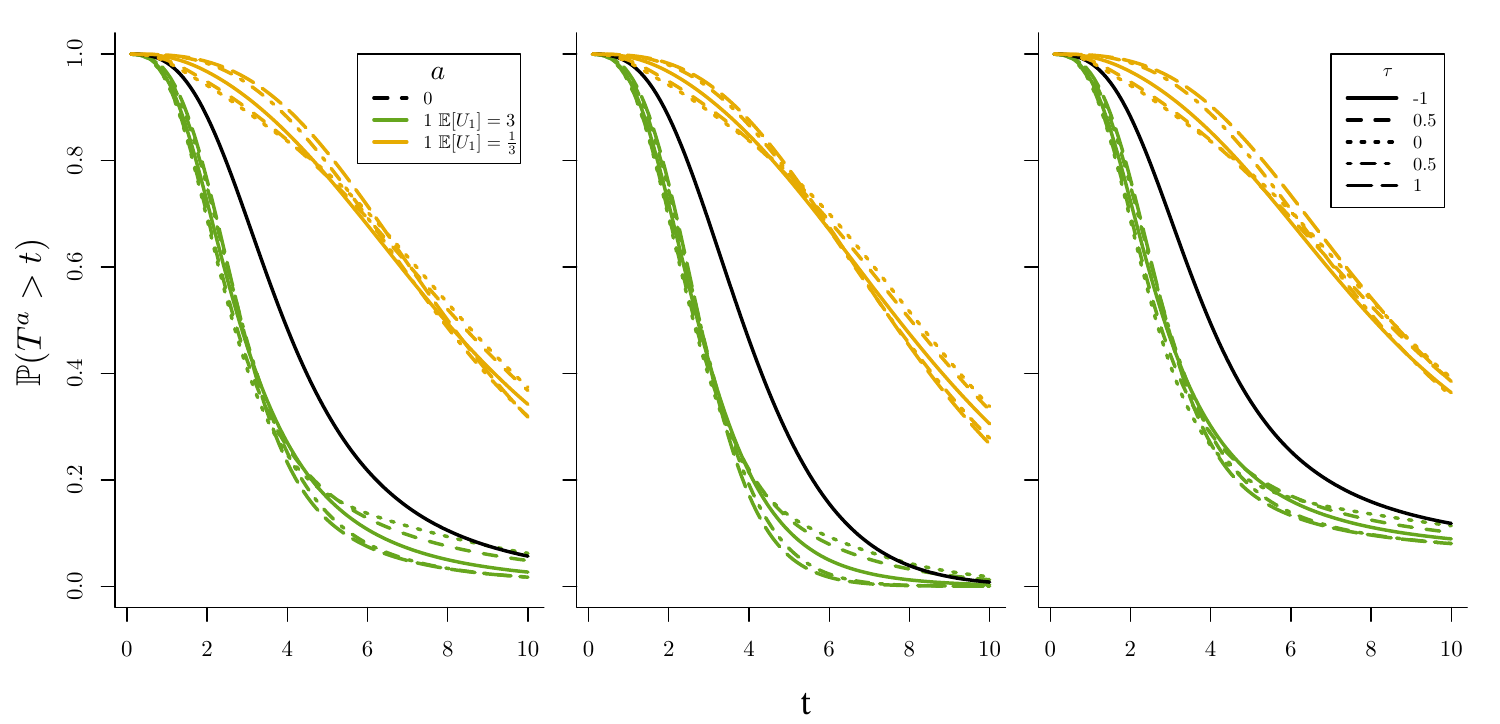} 
	\caption{Survival curves for $T^{0}$ (black), and $T^{1}$ (colored) when $\lambda^{a}_{i}(t)=U_{0i}(U_{1i})^{a} \tfrac{t^{2}}{20}$, for a unit-variance BHN distributed $U_{1}$ with $\mathbb{E}[U_{1}]=3$ (green) or $\mathbb{E}[U_{1}]=\tfrac{1}{3}$ (orange), when $U_{0}$ follows a gamma, inverse Gaussian or compound Poisson distribution (from left to right) with variance $1$ and where the joint distribution of $U_{0}$ and $U_{1}$ follows from a Gaussian copula with varying Kendall's $\tau$. The corresponding \mchr curves were presented in Figure \ref{CH5Fig5}. }\label{CH5Surv4}
\end{figure} 
\begin{figure}
	\centering
	\captionsetup{width=\textwidth}
	\includegraphics[width=0.75\textwidth]{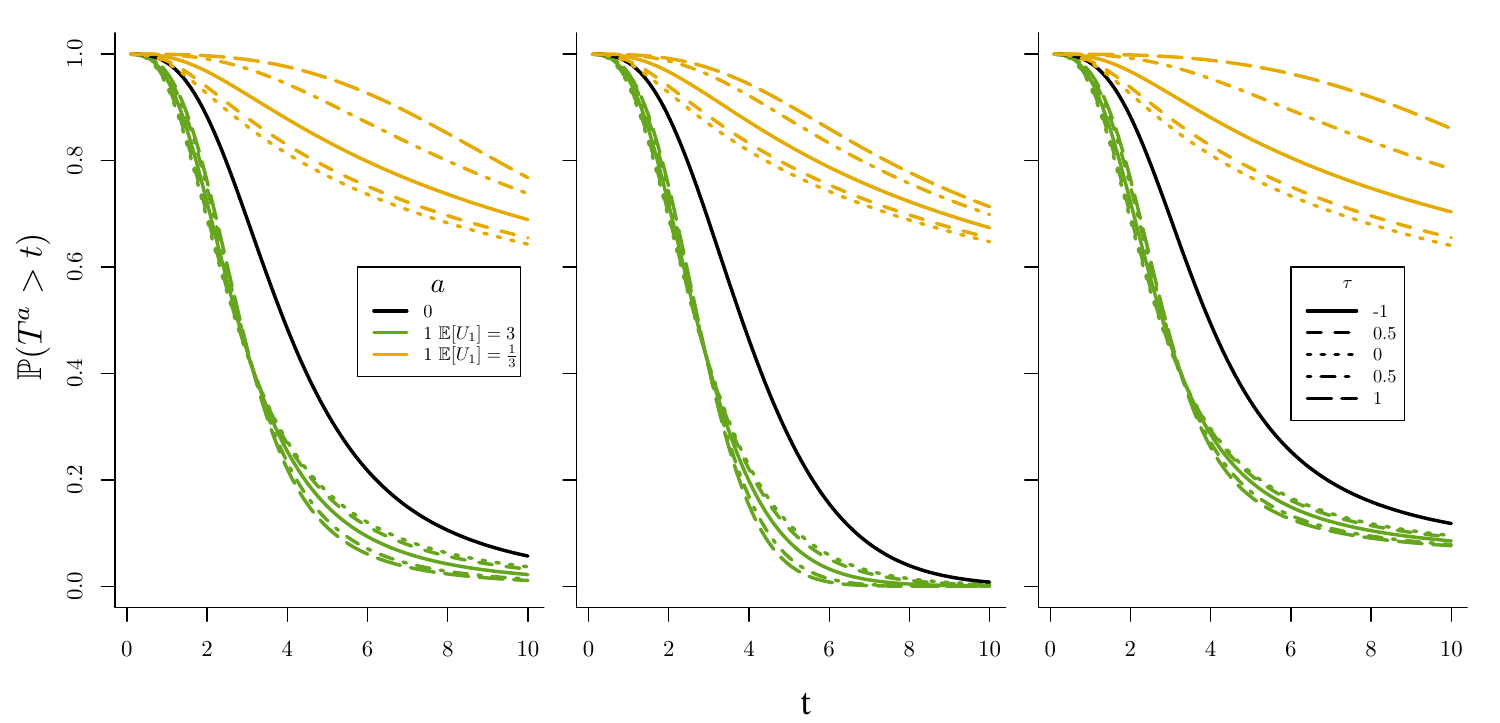} 
	\caption{Survival curves for $T^{0}$ (black), and $T^{1}$ (colored) when $\lambda_{i}^{a}(t)=U_{0i}(U_{1i})^{a} \tfrac{t^{2}}{20}$, for a unit-variance gamma distributed $U_{1}$ with $\mathbb{E}[U_{1}]=3$ (green) or $\mathbb{E}[U_{1}]=\tfrac{1}{3}$ (orange), when $U_{0}$ follows a gamma, inverse Gaussian or compound Poisson distribution (from left to right) with variance $1$. and where the joint distribution of $U_{0}$ and $U_{1}$ follows from a Gaussian copula with varying Kendall's $\tau$. The corresponding \mchr curves were presented in Figure \ref{CH5Fig5b}. }\label{CH5Surv5}
\end{figure}

	\end{document}